\documentclass[preprint,12pt]{elsarticle}

\usepackage{amssymb}
\usepackage{epsfig, graphicx, lscape, amsmath, amsbsy, amstext}
\usepackage{multirow}

\biboptions{authoryear}

\journal{Computational Statistics \& Data Analysis }

\def\b0{{\bf 0}}

\newcommand{\bl}{\begin{flushleft}}
\newcommand{\el}{\end{flushleft}}

\newcommand{\bc}{\begin{center}}
\newcommand{\ec}{\end{center}}

\begin{document}

\newtheorem{theorem}{Theorem}[section]
\newtheorem{prop}{Proposition}[section]
\newtheorem{corollary}{Corollary}[section]
\newtheorem{lemma}{Lemma}[section]

\begin{frontmatter}

\title{Regularized quantile regression under heterogeneous sparsity with application to quantitative genetic traits}


\author[ad1]{Qianchuan He\fnref{eq1}\corref{ca}}
\cortext[ca]{Corresponding author}
\ead{qhe@fhcrc.org}
\address[ad1]{Public Health Sciences Division, Fred Hutchinson Cancer Research Center, Seattle, WA 98109, U.S.A.}

\author[ad2]{Linglong Kong\fnref{eq1}}
\address[ad2]{Department of  Mathematical and Statistical Sciences, University of Alberta, Edmonton, AB Canada T6G 2G1}
\fntext[eq1]{Equally contributed.}

\author[ad3]{Yanhua Wang}
\address[ad3]{School of  Mathematics, Beijing Institute of Technology, Beijing, 100081, China}

\author[ad4]{Sijian Wang}
\address[ad4]{Department of Biostatistics \& Medical Informatics, The University of Wisconsin-Madison, Madison, WI 53706, U.S.A.}

\author[ad5]{Timothy A. Chan}
\address[ad5]{Human Oncology and Pathogenesis Program, Memorial Sloan-Kettering Cancer Center,  New York, NY 10065, U.S.A.}

\author[ad6]{Eric Holland}
\address[ad6]{Human Biology Division, Fred Hutchinson Cancer Research Center, Seattle, WA 98109, U.S.A.}

\begin{abstract}
Genetic studies often involve quantitative traits. Identifying genetic features that influence quantitative traits can help to uncover the etiology of diseases. 
Quantile regression method considers the conditional quantiles of the response variable, and is able to characterize the underlying regression structure in a more comprehensive manner. On the other hand, genetic studies often involve high-dimensional genomic features, and the underlying regression structure may be heterogeneous in terms of both effect sizes and sparsity. To account for the potential genetic heterogeneity, including the heterogeneous sparsity, a regularized quantile regression method is introduced.  The theoretical property of the proposed method is investigated, and its performance is examined through a series of simulation studies. A real dataset is analyzed to demonstrate the application of the proposed method.
\end{abstract}


\end{frontmatter}

\section{Introduction}

In many genetic studies, quantitative traits are collected for studying the associations between the traits and certain genomic features.
For example, body mass index, lipids and blood pressure have been investigated with respect to single nucleotide polymorphisms (SNPs)  \citep{avery2011phenomics}. With the rapid progress of high-throughput genome technology, new types of quantitative traits have emerged and attracted considerable research interest, such as the gene expression, DNA methylation, and protein quantification \citep{landmark2013genome}. The analysis of these quantitative traits yields new insight into biological processes and sheds light on the genetic basis of diseases.    


Typically, quantitative genetic traits are analyzed by least-square based methods, which seek to estimate the $E(Y|Z)$, where $Y$ is the trait and $Z$ is the set of covariates of interest. Quantile regression \citep{koenker1978regression} instead considers the conditional quantile function of $Y$ given $Z$, $Q_{\tau}(Y|Z)$, at a given $\tau \in (0,1)$.
When $\tau$ is fixed at 0.5, quantile regression is simply the median regression, which is well known to be more robust than the least-square estimation. In examining quantiles at different $\tau$, quantile regression provides a more complete picture of the underlying regression structure between $Y$ and $Z$.

Like least-square methods, traditional quantile regression methods only consider a handful of covariates. With the emergence of high-dimensional data, penalized quantile regression methods have been developed in recent years, and can be broadly classified into two classes. The first class seeks to harness the information shared among different quantiles to jointly estimate the regression coefficients. 
	\cite{jiang2014interquantile} proposed two novel methods, Fused Adaptive Lasso (FAL) and  Fused Adaptive Sup-norm estimator (FAS), for variable selection in interquantile regression. FAL combines the LASSO penalty \citep{tibshirani1996regression} with the fused lasso penalty \citep{tibshirani2005sparsity}, while FAS imposes the grouped sup-norm penalty and the fused lasso penalty. The fused lasso penalty delivers the effect of smoothing the regression slopes of adjacent quantiles, hence FAL and FAS can be used to identify common quantile slopes if such a smoothing property is desired. \cite{zou2008regularized} adopted an $F_{\infty}$ penalty,  which either eliminates or retains all the regression coefficients for a covariate at multiple quantiles. The method by \cite{jiang2013interquantile} seeks to shrink the differences among adjacent quantiles by resorting to the fused lasso penalty; however, this method does not perform variable selection at the covariate level, i.e., it does not remove any covariates from the model.

The second class of methods focuses on a single quantile at a time.
\cite{koenker2004quantile} imposed a LASSO penalty on the random effects in the mixed-effect quantile regression model; \cite{li2008l1} adopted the LASSO penalty; \cite{wu2009variable} explored the SCAD penalty \citep{fan2001variable}  and the adaptive LASSO penalty \citep{zou2006adaptive}, and proved the selection consistency and the normality of the proposed estimators (for a fixed dimension of the covariates).  \cite{wang2012quantile} investigated several penalties for quantile regression under the scenario of $p>n$, i.e., the dimension is larger than the sample size, and proved the selection consistency  of their proposed methods through a novel use of the subgradient theory. 
A recent approach by \cite{peng2014shrinkage} shares characteristics with both classes; its loss function targets a single $\tau$, while its penalty borrows information across different quantiles; their proposed penalty was shown to achieve more accurate results than the one that uses information from only a single quantile.


If the regression coefficients associated with a covariate are treated as a group, then some groups may be entirely zero and some other groups may be partially zero. Thus, sparsity can occur both at the group level and within the group level, and we call this type of sparsity as the heterogeneous sparsity. 
In this paper, we propose an approach that conducts joint variable selection and estimation for multiple quantiles under the situation that $p$ can diverge with $n$.  Our proposed method is able to achieve sparsity both at the group level and within the group level. 
We note that FAL can potentially yield sparsity at the two levels, but this approach has not been evaluated under the scenario where the dimension $p$ is high. 
To the best of our knowledge, this is the first paper that explicitly investigates the heterogeneous sparsity for quantile regression.
We show that our method tends to be more effective than the compared methods in handling heterogeneous sparsity when the dimension is high. We also provide theoretical justification for the proposed method. The paper is organized as the follows. In Section 2, we describe the proposed method and the implementation details. In Section 3, we prove the theoretical properties of the proposed method. In Section 4, we show the results of simulation studies regarding several related methods, and in Section 5, we present an example of real data analysis for the proposed method.

\section{Method}
\subsection{Data and Model}

Let $Z$ represent the vector consisting of $p$ covariates, such as SNPs or genes. Let $\gamma_\tau$ be the $p$-dimension coefficient vector at the $\tau$th quantile. Let $Y$ be the random variable that denotes the phenotype we are interested in, such as quantitative traits in genetic studies. For a given $\tau \in (0,1)$, the linear quantile regression model is known as
\begin{equation*}
	Q_{\tau}(Y|Z)= \gamma_0 + Z^{\rm T}\gamma_\tau ,
\end{equation*}
where $\gamma_0$ is the intercept, and $Q_{\tau}(Y|Z)$ is the $\tau$-th conditional quantile of $Y$ given $Z$, that is, $P_{Y|Z}(Y\leq Q_{\tau}(Y|Z))=\tau$.

The dimension
$p$ can be potentially very high in genomic studies, but typically it is assumed that only a limited number of genomic features contribute to the phenotype. 
For this reason, one needs to find a sparse estimation of  $\gamma_\tau$ to identify those important genomic features.
On the other hand, we also wish to consider multiple quantile levels simultaneously so that information shared among different quantile levels can be utilized. To this end, we propose the following model for the joint estimation of the regression coefficients for multiple quantiles.
Given $M$ quantile levels, $0<\tau_1<\cdots<\tau_M<1$, our linear quantile regression model is defined as, for $\tau_m$ $(m=1,\cdots,M)$, 
\begin{align}\label{0-mod-1}
	Q_{\tau_m}(Y|Z)=\gamma_{m0} + Z^{\rm T}\gamma_{\tau_m}.
\end{align}
where $\gamma_{m0}$ is the intercept, and $\gamma_{\tau_m} $ is the $p$-dimension coefficient vector. For ease of notation, we write $\gamma_m=\gamma_{\tau_m}$.
For the above model, we further define $\gamma \equiv (\gamma_1^{\rm T}, \cdots,\gamma_M^{\rm T})^{\rm T}$ 
with $\gamma_m= (\gamma_{m1}, \cdots,\gamma_{mp})^{\rm T}$,
and intercept parameter $\gamma_0 \equiv (\gamma_{10},\cdots,\gamma_{M0})^{\rm T}$. 

We now focus on the sample version of the model \eqref{0-mod-1}. Let $\{(Y_i,Z_i^{\rm T})^{\rm T} \}_{i=1}^n$ be an i.i.d. random sample of size $n$ from population $(Y,Z^{\rm T})^{\rm T}$, where $Z_i=(Z_{i1},Z_{i2},\cdots,Z_{ip})^{\rm T}$.
 The sample quantile loss function is defined as
\begin{align*}
	 Q_n(\gamma_0, \gamma)=\sum_{m=1}^M\sum_{i=1}^{n}\,\rho_m(Y_i-Z_i^{\rm T}\gamma_m-\gamma_{m0})
\end{align*}
where $\rho_m(u)=u(\tau_m-I(u<0))$ is the quantile check loss function with $I(\cdot)$ being the indicator function.
To introduce sparsity to the model, we add to the loss function a penalty function 
\begin{align*}
	P_n(\gamma) = n\lambda_n\sum_{j=1}^p\left(\sum_{m=1}^M\,\omega_{mj}|\gamma_{mj}|\right)^{\frac{1}{2}},
\end{align*}
where $\lambda_n$ is the tuning parameter, and $\omega_n=(\omega_{mj}:1\leq m\leq M, 1\leq j\leq p)$ is the weight vector whose component   $\omega_{mj}>0$ is the weight of parameter $\gamma_{mj}$. Note that the penalty is a nonconvex function. It essentially divides the regression coefficients into $p$ groups, and each group consists of $M$ parameters associated with the $j$th covariate. The motivation is that, while each
quantile may have its own set of regression parameters, we wish to borrow strengths from each
quantiles to select covariates that are important across all quantiles
as well as those covariates that are important to only some of the
quantiles. This type of nonconvex penalty has been considered in the Cox regression model and other settings; see \cite{wang2009hierarchically} for an example, but to the best of our knowledge has not been studied in the quantile regression model. 
We can choose
$\omega_{m j}= \tilde\gamma_{mj}^{-1}$, where $\tilde\gamma_{m
j}$ is some consistent estimate for $\gamma_{mj}$. For example, we may use the estimates from the unpenalized quantile regression conducted at each individual quantile level. When $p<n$ but is fixed, the consistency of the unpenalized estimates has been proved by \citet{koenker1978regression}. When $p<n$ but is diverging with $n$, the estimates from unpenalized quantile regression are consistent by adapting to Lemma A.1 of \citet{wang2012quantile}.  Thus, our objective function is defined as
$$
\sum_{m=1}^M\sum_{i=1}^{n}\,\rho_m(Y_i-Z_i^{\rm T}\gamma_m-\gamma_{m0})
	 + n\lambda_n\sum_{j=1}^p\left(\sum_{m=1}^M \omega_{mj}|\gamma_{mj}|\right)^{\frac{1}{2}}.
$$

For the sake of convenience, we define
$	\theta =(\theta_1^{\rm T},\cdots,\theta_M^{\rm T})^{\rm T} \quad \text{with}\quad \theta_m = (\gamma_{m0},\,\gamma_m^{\rm T})^{\rm T},\quad m=1,\cdots,M,$ 
and the corresponding parameter space by $\Theta_n\subset \mathbb{R}^{M(p+1)}$.
Further define $U_i =(1,Z_i^{\rm T})^{\rm T}$, $i=1,\cdots,n$.
Then, $Q_n(\gamma_0, \gamma)$ can be written as $Q_n(\theta)$. Emphasizing that $\gamma$ is a subset of $\theta$, we can write the objective function as
\begin{align}\label{0-ln}
	L_n(\theta) & \equiv Q_n(\theta) + P_n(\gamma)
	 =\sum_{m=1}^M\sum_{i=1}^{n}\,\rho_m(Y_i-U_i^{\rm T}\theta_m)
	 + n\lambda_n\sum_{j=1}^p\left(\sum_{m=1}^M \omega_{mj}|\gamma_{mj}|\right)^{\frac{1}{2}}.
\end{align}

Let $\hat{\theta}$ be a local minimizer of $L_n(\theta)$ in \eqref{0-ln} for $\theta\in\Theta_n$. Because the heterogeneity of sparsity is explicitly taken into account in this model, we name our proposed method as Heterogeneous Quantile Regression (Het-QR). 
Our model can be modified to accommodate different weights for the losses at different quantiles. That is, $Q_n(\theta)$ may take the form $\sum_{m=1}^M \pi_m \sum_{i=1}^{n}\,\rho_m(Y_i-U_i^{\rm T}\theta_m)$, where $\pi_m$ is the weight for the $m$th quantile. Some examples on the choice of weight $\pi_m$ can be found in \cite{koenker2004quantile} and \cite{zhao2014efficient}.   

%

\subsection{Implementation}
We design the following algorithm to implement the proposed method.
First, we show in the Appendix that the objective function can be transformed into
$${\rm argmin}_{  \theta,  \xi}   \sum_{m=1}^M\sum_{i=1}^{n}\,\rho_m(Y_i-U_i^{\rm T}\theta_m)+
\lambda_1\sum_{j=1}^p \xi_j +
 \sum_{j=1}^p \xi_j^{-1} \left( \sum_{m=1}^M \omega_{mj}
|\gamma_{mj}| \right), $$
where $\xi=(\xi_1,\ldots,\xi_p)$ are newly introduced nonnegative parameters. 
Then, the new objective function can be solved by the following iterative algorithm:\\
{\bf Step 1:} We first fix $\theta$ to solve $\xi_j, j=1,...,p$. To this end,
$\xi_j$ has a closed-form solution. That is,
$\hat\xi_j=\left(\sum_{m=1}^M \omega_{mj} |\gamma_{mj}|
\right)^{1/2} \lambda_1^{-1/2}, j=1,...,p$.
 \\
{\bf Step 2:} We fix $\xi_j, j=1,...,p$, to solve
$\theta$. That is, we aim to solve 
$${\rm argmin}_{  \theta}   \sum_{m=1}^M\sum_{i=1}^{n}\,\rho_m(Y_i-U_i^{\rm T}\theta_m)
+
 \sum_{j=1}^p   \sum_{m=1}^M \xi_j^{-1} \omega_{mj}
|\gamma_{mj}|.$$ 
We can formulate this objective function as a linear program and derive its dual form (see Appendix for details),  then
the optimization can be conducted by recoursing to
 existing linear programming packages; we utilize the Quantreg R package \citep{koenker2015quantile}.
\\
{\bf Step 3:} Iterate step 1 and step 2 until convergence. Due to the nonconvexity of the penalty function, the estimate is a local minimizer.

\section{Theoretical Properties}

Now we investigate the asymptotic properties of the proposed method. FAL and FAS considered $p$ to be fixed. We study the situation where $p$ can diverge with $n$.
Let the true value of $\theta$ be $\theta^*$, 
where the corresponding true values of $\gamma_{m0}, \gamma_{mj}, \gamma$ are $\gamma_{m0}^*, \gamma_{mj}^*, \gamma^*$, respectively. 
Let the number of nonzero elements in  $\gamma^*$ be $s$.
To emphasize that $s$ and $p$ can go to infinity, we use $s_n$ and $p_n$ when necessary. For Theorems 1 and 2 (to be shown), $p_n$ is at the order lower than $O(n^{1/2})$; for Theorem 3, $p_n$ is at the order lower than $O(n^{1/6})$.

We define some index sets to be used in our theorems. 
Let $\mathcal{N}=\{(m,j):1\leq m\leq M,1\leq j\leq p_n\}$.
For the true parameters, define the oracle index set $I=\{(m,j)\in \mathcal{N}:\gamma_{mj}^* \neq 0\}$
and its complementary set $II =\{(m,j) \in \mathcal{N}:\gamma_{mj}^* = 0\}$.
Assume that $I$ has cardinality $|I|=s_n$.


We define some notations used for our theorems. Define $d_{nI}=\max_{(m,j)\in I}\omega_{mj}$ and 
$d_{nII}=\min_{(m,j)\in II} \omega_{mj} \bigl( \max_{(m,j)\in\mathcal{N}}\omega_{mj} \bigr)^{-\frac{1}{2}}$.
Define $\theta_I^*$ and $\hat{\theta}_I$ as the subvectors of the vectors $\theta^*$ and $\hat{\theta}$ corresponding to the oracle index set $I$, respectively.
For every fixed $1\leq m\leq M$, define the index set $I_m = \{1\leq j\leq p_n:\gamma_{mj}^* \neq 0\}$. Let
\begin{align}
	\mathbf{\Sigma}_n &= ( \mathbf{\Sigma}_{lm} )_{M\times M} \quad\text{with}\quad 
		\mathbf{\Sigma}_{lm} = (\min(\tau_m,\tau_l) - \tau_m\tau_l)\, E(U_{ilI}U_{imI}^{\rm T}),\label{4-sigma}\\
	\mathbf{B}_{nI} &= Diag(B_{1},B_{2},\cdots,B_{M}) \quad \text{with}\quad 
		B_{m}=\sum_{i=1}^n\,f(U_i^{\rm T}\theta_m^*|\,Z_i)\,U_{imI}U_{imI}^{\rm T} \label{4-bn},
\end{align}
where $U_{imI} =(1,Z_{imI}^{\rm T})^{\rm T}$ with $Z_{imI}$ being the subvector of $Z_i$ corresponding to index set $I_m$.


Let $F(y|z)$ and $f(y|z)$ be the conditional distribution function and the conditional density function of $Y$ given $Z=z$, respectively. 
For any proper square matrix $A$, let $\lambda_{min}(A)$ and $\lambda_{max}(A)$ denote the minimum and maximum eigenvalue of $A$, respectively.

Before stating the main theorems, we need the following regularity conditions labeled by $\mathcal{L}$:
\renewcommand{\labelenumi}{(L\arabic{enumi})}
\begin{enumerate}
    \item The conditional density $f(y|z)$  has first order derivative $f'(y|z)$ with respect to $y$;
    And $f(y|z)$ and $f'(y|z)$ are uniformly bounded away from $0$ and $\infty$ on the support set of $Y$ and the support set of $Z$;

    \item For randon sample $Z_i=(Z_{i1},Z_{i2},\cdots,Z_{ip})^{\rm T}$, $1\leq i\leq n$, there exists a postive constant $C_1$ 
    such that $\max_{1\leq i\leq n,1\leq j\leq p}|Z_{ij}|\leq C_1$;

    \item For $U_i=(1,Z_i^{\rm T})^{\rm T}$, $i=1,\cdots,n$, let $\mathbf{S}_n=\sum_{i=1}^n\,U_iU_i^{\rm T}$. There exist positive constants $C_2<C_3$ such that 
    $C_2\leq \lambda_{min}(n^{-1}\mathbf{S}_n) \leq \lambda_{max}(n^{-1}\mathbf{S}_n)\leq C_3$;

    \item The dimension $s_n$ satisfies that $s_n= a_0n^{\alpha_0}$, and the dimension $p_n$ satisfies that $p_n=a_1n^{\alpha_1}$, where $0<\alpha_0 <\alpha_1<\frac{1}{2}$, and $a_0$ and $a_1$ are two positive constants;
    
    \item The matrix $\mathbf{B}_n$ given in $\eqref{0-anbn}$ (see Appendix) satisfies that
    $C_4\leq \lambda_{min}(n^{-1}\mathbf{B}_n)\leq \lambda_{max}(n^{-1}\mathbf{B}_n)\leq C_5$, where $C_4$ and $C_5$ are positive constants. 
 
    \item  The matrix $\mathbf{\Sigma}_n$ satisfies that $\lambda_{min}(\mathbf{\Sigma}_n) \geq C_6$, where $C_6$ is a positive constant.
 
\end{enumerate}


Conditions (L1)-(L3) and (L5)-(L6) are seen in typical theoretical investigation of quantile regression. Condition (L4) specifies the magnitude of $s_n$ and $p_n$ with respect to the sample size.   Under the aforementioned regularity conditions, we present the following three theorems. The proof is relegated to the Appendix. Define $\theta_{II}^*$ and $\hat{\theta}_{II}$ as the subvectors of the vectors $\theta^*$ and $\hat{\theta}$ corresponding to the index set $II$, respectively. Clearly, $\theta_{II}^*=0$. Due to the nonconvexity of the penalty function, all the following theorems and their proofs (in the Appendix) are regarding to a local minimizer of the objective function.

{\sc Theorem 1.}
{\it Under conditions (L1)-(L5),
if $\lambda_nd_{nI}^\frac{1}{2}=o(s_n^{-\frac{3}{4}}p_n^{\frac{3}{4}}n^{-\frac{3}{4}})$, then
the estimator $\hat{\theta}$ of $\theta^*$ exists, is a local minimizer, and satisfies the estimation consistency that
$\|\hat{\theta}-\theta^*\|_2 = O_{p}(n^{-\frac{1}{2}}p_n^{\frac{1}{2}})$.}

Theorem 1 shows that the proposed method is consistent in parameter estimation. The convergence rate $O_{p}(n^{-\frac{1}{2}}p_n^{\frac{1}{2}})$ is typical for the settings where $p$ diverges with $n$. 

{\sc Theorem 2.}
{\it Under conditions (L1)-(L5),
if $\lambda_n d_{nI}^\frac{1}{2}=o(n^{-\frac{3}{4}}s_n^{-\frac{3}{4}}p_n^{\frac{3}{4}})$ 
and $n^{-\frac{1}{2}}p_n^{\frac{1}{2}}=o(\lambda_n\,d_{nII})$, then $P(\hat{\theta}_{II}=0)\to 1$. }

Theorem 2 indicates that our method can distinguish the truly zero coefficients from the nonzero coefficients with probability tending to 1. It can be seen that the penalty weight $d_{nII}$ plays a critical role  in the property of selection consistency.

{\sc Theorem 3.} 
{\it Under conditions (L1)-(L3) and (L5)-(L6), if $\lambda_n d_{nI}^{\frac{1}{2}}=O(n^{-1}p_n^{\frac{1}{2}})$, and $n^{-\frac{1}{2}}p_n^{\frac{1}{2}}=o(\lambda_n\,d_{nII})$, and the powers of $s_n$ and $p_n$ in condition $(L4)$ satisfy $0<\alpha_0<\alpha_1<\frac{1}{6}$, then for any unit vector $b\in\mathbb{R}^{M+s_n}$ 
 we have
\begin{align*}
	(nb^{\rm T} \mathbf{\Sigma}_n b)^{-\frac{1}{2}}\,b^{\rm T}\mathbf{B}_{nI}\bigl(\hat{\theta}_I-\theta_I^* \bigr)\to N(0,1).
\end{align*}
}
Theorem 3 suggests that the estimated nonzero coefficients have the asymptotic normality.
Heuristically, for given $n, \lambda_n$  and $\omega_{mj}$, the considered penalty has its slope tending to infinity when $\gamma_{mj}$ goes to 0, thus the penalty tends to dominate small $\gamma_{mj}$. On the other hand, when $\lambda_n$ is sufficiently small, the penalty has little impact on the estimation of relatively large $\gamma_{mj}$. These properties, in combination with proper choice of the tuning parameter, play major roles in the oracle property of the proposed estimator. 
The oracle property for coefficients within a group is mainly due to the penalty weights, which put large penalty on small coefficients (and small penalty on large coefficients).

%
%
%

\section{Simulation Studies}

We conduct simulation studies to evaluate the proposed method along with the following methods: the QR method, which applies  quantile regression to each individual quantile level without any variable selection; the QR-LASSO method, which adopts the $L_1$-penalized quantile regression to each quantile level; the QR-aLASSO method, which imposes the adaptive LASSO penalty to each quantile level \citep{wu2009variable}; the FAL and the FAS method \citep{jiang2014interquantile}. Both FAL and FAS contain a fused-LASSO type of penalty, which encourages the equality of the regression coefficients among different quantiles. FAL allows within-group sparsity,   while FAS generates sparsity only at the group level. For Het-QR, we set the penalty weight $\omega_{m j}$ to be the inverse of the estimate from the unpenalized quantile regression (unless specified otherwise).

We first consider a model where important covariates have nonzero regression coefficients across all (or almost all) quantiles. We simulate 6 independent covariates, each of which follows the uniform(0,1) distribution. Then we simulate the trait as
$$Y=1.0+ \beta_1 Z_1 + \beta_2 Z_2 + \beta_6 Z_6 + \kappa Z_6 \epsilon,$$   
where  $\beta_1=1, \beta_2=1, \beta_6=2$, $\kappa=2$ and $\epsilon \sim N(0,1)$.  Under this set up, $Z_1$ and $Z_2$ have constant regression coefficients across all quantiles, while $Z_6$'s regression coefficient is determined by $2+ 2 \times\Phi^{-1}(\tau)$, which varies among different quantiles. 
That is, the $\tau$th quantile of $Y$ given $Z_1, Z_2$ and $Z_6$ is 
$$Q_{\tau}\left(Y|Z_1,Z_2,Z_6\right)=1.0+ Z_1 + Z_2 + \left(2+ 2 \times\Phi^{-1}(\tau)\right)Z_6.$$  
 This model is in line with the model considered by \cite{jiang2014interquantile}. All the other 3 covariates, $Z_3, Z_4$ and $Z_5$, have no contribution to $Y$. The sample size $n$ is set to 500. To select the tuning parameter, we follow the lines of \cite{mazumder2011sparsenet} and \cite{wang2012quantile} to generate another dataset with sample size of $10n$, and then pick the tuning parameter at which the check loss function  is minimized. The total number of simulations for each experiment is $100$. 

We consider various criteria to evaluate the performance of the compared methods, such as the model size and 
the parameter estimation error (PEE).  The model size refers to the number of estimated non-zero coefficients among the $M$ quantile levels.  The PEE is calculated by $ \sum_{m=1}^M \sum_{j=1}^p |\hat \gamma_{mj} -\gamma^*_{mj}|/M$.
To evaluate the prediction error, we simulate an independent dataset, $(Y_{pred}, Z_{pred})$, with sample size of 100n, and then calculate the F-measure (FM) \citep{gasso2009recovering}, the quantile prediction error (QPE) and the prediction error (PE). The FM is equal to $2\times S_a/M_a$, where $S_a$ is the number of truly nonzero slopes being captured, and $M_a$ is the sum of the estimated model size and the true model size.
The QPE is defined as the sample version of the $\sum_{m=1}^M (Q_{\tau_m}(Y_{pred}|Z_{pred}) - Z_{pred}^{\rm T}\hat\gamma_{\tau_m}-\hat\gamma_{m0})^2/M$, averaged across all the subjects. The PE is defined as $Q_n(\hat\theta)/n$, i.e., the check loss averaged across all considered quantiles for all the test samples,  evaluated on $(Y_{pred}, Z_{pred})$. 

For the purpose of illustration, we consider three quantiles, $\tau=0.25, 0.5, 0.75$. The results are shown in the upper panel of Table 1. It can be seen that when $p$ is 6, FAL has the lowest parameter estimation error and  FAS has the lowest PE, though the difference between these two methods and the other compared methods is generally quite small. Next, we increase $p$ to 100 to evaluate the the methods under higher dimension. As shown in the lower panel of Table 1, when $p$ is equal to 100, both FAL and FAS have deteriorated performance; for instance, their model sizes tend to be twice (or more) as the true model size and their PEE and PE are higher than Het-QR. This experiment shows that the performance of FAL and FAS is suboptimal when the dimensionality grows large; one potential explanation is that the penalties of FAL and FAS may overemphasize the interquantile shrinkage, which make them less efficient when many noise covariates are present. Further research is merited. As to computation, we did not observe non-convergence for Het-QR in our experiments.

\begin{table}
Table 1. Comparison of Het-QR and other methods in the absence of within-group sparsity (standard error of the sample mean shown in the parenthesis 
)\small 

\begin{center}

 {\begin{tabular}{cccccccc}
\hline\hline
  Method  & Model-size & FM (\%)  & $PEE\times 100$ &  $QPE\times 10^3$ 
   & $PE \times 10^3$&  
\\ \hline
 & \multicolumn{5}{c} {$p=6$ }
\\
QR           &18            &-        &53.3(1.5)  &10.2(0.6)  &1041.2(0.6)   & \\
QR-LASSO     &15.0(0.2)     &75(0.7)   &36.1(1.2)  &6.7(0.5)   &1038.4(0.6)   &\\ 
QR-aLASSO    &10.9(0.2)     &91(0.7)   &25.6(0.9)  &5.8(0.4)   &1037.2(0.5)   &\\ 
FAL          &11.1(0.2)     &91(0.9)   &25.3(1.0)  &6.2(0.5)   &1037.2(0.5)   &\\      
FAS          &12.1(0.2)     &86(0.9)   &26.4(1.1)  &5.9(0.4)   &1037.1(0.5)   &\\
Het-QR       &9.6(0.1)      &97(0.6)   &26.0(0.9)  &6.5(0.5)   &1037.5(0.5)   &\\ 
 & \multicolumn{5}{c} {$p=100$ } 
\\
QR           &300           &-           &1556.4(12.7) &325.0(5.1) &1242.1(2.7) &\\
QR-LASSO     &47.3(1.2)     &33(0.7)   &120.5(3.5)   &23.4(1.2)  &1052.4(0.9) &\\ 
QR-aLASSO    &16.8(0.4)     &72(1.1)   &46.9(1.9)    &10.5(0.8)  &1042.1(0.7) &\\ 
FAL          &17.7(0.6)     &70(1.4)   &41.9(1.8)    &10.2(0.8)  &1041.0(0.6) &\\      
FAS          &23.7(0.7)     &58(1.2)   &58.0(2.4)    &13.9(0.9)  &1043.4(0.7) &\\
Het-QR       &9.3(0.1)      &99(0.4)   &29.5(1.2)    &8.4(0.6)   &1039.9(0.6) &\\ 
 \hline
\end{tabular}}
\\
{}
\end{center}
\end{table}

Next, we systematically evaluate the situation where within-group sparsity exists. To introduce correlations into covariates, we simulate 20 blocks of covariates, each block containing 5 correlated covariates. For each block, we first simulate a multivariate
normal distribution with mean being the unit vector and covariance
matrix following either the compound symmetry or the
auto-regressive correlation structure with correlation coefficient $\rho=0.5$; next, we take the absolute value of the simulated random normal variables as the covariates $Z$. The total number of covariates is 100. 
We specify the conditional quantile regression coefficient function $\gamma(\tau)$  as follows.
For $\tau \in (0, 0.3]$, the first 8 regression slopes for $Z$ are $(0.5, 0, 0,0,0, 0.6, 0, 0)$; for  $\tau \in (0.3, 0.7]$, the first 8 regression slopes are $(0.5, 0, 0, 0, 0, 0.6, 0, 0.7)$; for $\tau \in (0.7, 1.0)$, the corresponding slopes are $(0.6, 0,0,0,0, 0.7, 0, 0.7)$. All other regression slopes are 0. Thus, the first and the sixth covariates are active among all quantiles, while the eighth covariate is active only for the last two quantile levels. To generate $Y$, we first simulate a random number $\tau \in \mbox{Uniform}(0,1)$, and then determine the $\gamma(\tau)$ based on $\tau$; subsequently, we obtain $$Y= 1.0+ Z^{\rm T} \gamma(\tau) + F^{-1} (\tau),$$ where $F^{-1}$ is the inverse cumulative function of some distribution $F$. That is, the $\tau$th quantile of $Y$ given $Z$ is 
\[
Q_{\tau}\left(Y|Z\right) = 
  \begin{cases} 
1.0+ F^{-1} (\tau) +0.5Z_1+0.6Z_6 & \text{if } 0<\tau \leq 0.3 \\
1.0+ F^{-1} (\tau) +0.5Z_1+0.6Z_6 +0.7 Z_8& \text{if } 0.3<\tau \leq 0.7 \\
1.0+ F^{-1} (\tau) +0.6Z_1+0.7Z_6 +0.7 Z_8& \text{if } 0.7<\tau <1 
  \end{cases}.
\]
We explore different distributions for $F$: the standard normal distribution, the $T$-distribution with degrees of freedom equal to 3 ($T_3$), and the exponential distribution with shape parameter equal to 1. 

We first consider the normal distribution for $F$.  The results are shown in Table 2. Because no variable selection is conducted, QR has much larger PEE, QPE, and PE than the other methods; for example, the PEE and QPE of QR are more than 10 times higher than the compared methods. 
QR-LASSO, QR-aLASSO, FAL, and FAS have more tamed model sizes, but still contain a number of noise features. Het-QR yields a model that is closer to the true model, in which the three considered quantiles contains 2, 2, and 3 nonzero slopes, respectively. 
 Het-QR also appears to have the highest FM, and lowest errors for parameter estimation and prediction. Next, we consider the distribution to be the $T_3$ (Table 3) and the exponential distribution (Table 4), and the results show a similar pattern. These experiments indicate that  Het-QR can handle higher dimension as well as the heterogeneous sparsity better than the other methods.

We finally consider the situation where $p>n$. While theoretical development is still needed for this setting, our experiment is to evaluate the practical performance of the proposed approach.  We let $n=500$ and $p=600$. For $\tau \in (0, 0.3]$, the first 8 regression slopes for $Z$ are $(0.6, 0, 0,0,0,  0.6, 0, 0)$; for  $\tau \in (0.3, 0.7]$, the first 8 regression slopes are $(0.6, 0, 0.8, 0, 0,   0.7, 0, 0.8)$; for $\tau \in (0.7, 1.0)$, the corresponding slopes are $(0.8, 0,0.8,0,0,   0.8, 0, 1.0)$. In this scenario, $Z_3$ and $Z_8$ have zero coefficients for the first quantile, but nonzero coefficients for the other two quantiles. That is, the $\tau$th quantile of $Y$ given $Z$ is 
\[
Q_{\tau}\left(Y|Z\right) = 
  \begin{cases} 
1.0+ F^{-1} (\tau) +0.6Z_1+0.6Z_6 & \text{if } 0<\tau \leq 0.3 \\
1.0+ F^{-1} (\tau) +0.6Z_1+0.8Z_3 +0.7 Z_6+0.8Z_8& \text{if } 0.3<\tau \leq 0.7 \\
1.0+ F^{-1} (\tau) +0.8Z_1+0.8Z_3 +0.8 Z_6+1.0Z_8& \text{if } 0.7<\tau <1 
  \end{cases}.
\]
We omit QR, FAL and FAS  because they are not designed to handle the setting of `$p>n$'. QR-LASSO can be directly applied to the data that have higher dimension than sample size.
For  QR-aLASSO, we derive the penalty weights using the estimates from QR-LASSO.  For Het-QR, we first run Het-QR with the penalty weights equal to 1 to obtain the initial estimators for $\gamma_m^*, m=1,\ldots, M$, and then use the inverse of the initial estimators as the penalty weights; finally, we run  Het-QR to obtain the $\hat \theta$. The results are shown in Table 5.  Het-QR tends to yield a smaller model  than the compared methods and have better performance in estimating the regression coefficients as well as in prediction. 



\begin{table}
Table 2. Comparison of Het-QR and other methods for $p=100$ under the normal distribution (standard error of the sample mean shown in the parenthesis
)\small

\begin{center}

 {\begin{tabular}{cccccccc}
\hline\hline
  Method  & Model-size  & FM (\%)& $PEE\times 100$ &  $QPE\times 10^3$ 
   & $PE \times 10^3$& 
\\ \hline
 & \multicolumn{5}{c} {Correlation structure: Auto-regressive }
\\
QR           &300        & -         &1189.3(7.2)  &791.6(9.1)  &1606.5(3.1)   &\\
QR-LASSO     &41.1(1.3)  &34(0.9)  &108.7(3.1)   &78.5(3.1)   &1353.1(1.2)   &\\ 
QR-aLASSO    &17.8(0.5)  &63(1.2)  &60.4(2.5)    &49.2(3.2)   &1341.5(1.2)   &\\ 
FAL          &20.3(0.8)  &60(1.6)  &59.5(2.6)    &49.9(3.5)   &1340.4(1.2)   &\\      
FAS          &24.8(1.1)  &53(1.5)  &85.1(2.6)    &86.4(3.1)   &1351.2(1.1)   &\\
Het-QR       &8.6(0.1)   &96(0.7)  &34.7(1.6)    &32.6(2.8)   &1334.9(1.1)   &\\ 
 & \multicolumn{5}{c} {Correlation structure: Compound symmetry } 
\\
QR           &300        & -         &1218.1(8.4)  &792.7(9.5)  &1606.1(3.4)   &\\
QR-LASSO     &40.1(1.1)  &35.0(0.9)  &105.9(3.1)   &74.8(3.0)   &1351.9(1.3)   &\\ 
QR-aLASSO    &18.9(0.6)  &61(1.3)  &61.2(2.7)    &48.0(3.1)   &1341.6(1.3)   &\\ 
FAL          &19.9(0.8)  &61(1.4)  &61.0(2.9)    &54.4(4.0)   &1341.2(1.3)   &\\      
FAS          &24.1(1.1)  &55(1.6)  &83.8(2.8)    &89.3(3.2)   &1351.3(1.2)   &\\
Het-QR       &8.9(0.2)   &94(0.8)  &34.6(1.9)    &31.4(2.9)   &1334.7(1.2)   &\\ 
 \hline
\end{tabular}}
\\
{}
\end{center}
\end{table}

\begin{table}
Table 3. Comparison of the Het-QR and other methods for $p=100$ under the $T_3$ distribution (standard error of the sample mean shown in the parenthesis
)\small

\begin{center}

 {\begin{tabular}{cccccccc}
\hline\hline
  Method  & Model-size  & FM (\%)& $PEE\times 100$ &  $QPE\times 10^3$ 
   & $PE \times 10^3$& 
\\ \hline
 & \multicolumn{5}{c} {Correlation structure: Auto-regressive }
\\
QR           &300        & -         &1452.1(10.2) &1149.1(15.2)  &2114.3(4.5)   &\\
QR-LASSO     &39.4(1.3)  &35(0.9)  &123.6(3.3)   &103.7(3.7)    &1796.8(1.4)   &\\ 
QR-aLASSO    &19.7(0.6)  &58(1.3)  &81.0(3.0)    &77.3(4.2)     &1787.7(1.6)   &\\ 
FAL          &20.6(0.7)  &58(1.3)  &76.4(3.1)    &74.8(4.6)    &1786.2(1.6)   &\\      
FAS          &24.2(0.9)  &53(1.3)  &99.6(3.4)    &107.8(4.2)   &1795.6(1.6)   &\\
Het-QR       &8.9(0.2)   &92(1.0)  &47.3(2.3)    &53.5(4.0)    &1779.4(1.5)   &\\ 
 & \multicolumn{5}{c} {Correlation structure: Compound symmetry } 
\\
QR           &300        & -         &1487.4(11.2)  &1156.3(15.6)  &2115.2(4.7)   &\\
QR-LASSO     &38.7(1.0)  &35(1.2)  &120.5(3.2)   &99.4(3.4)   &1795.5(1.6)   &\\ 
QR-aLASSO    &20.3(0.7)  &56(1.3)  &81.6(3.3)    &76.0(4.1)   &1787.6(1.8)   &\\ 
FAL          &22.5(0.9)  &56(1.5)  &82.5(3.9)    &78.1(5.1)   &1786.6(1.8)   &\\      
FAS          &25.6(1.0)  &51(1.3)  &99.1(3.5)    &107.0(4.1)  &1794.7(1.7)   &\\
Het-QR       &9.1(0.3)   &91(1.1)  &48.6(2.9)    &55.4(4.7)   &1780.1(1.8)   &\\ 
 \hline
\end{tabular}}
\\
{}
\end{center}
\end{table}

\begin{table}
Table 4. Comparison of Het-QR and other methods for $p=100$ under the exponential distribution (standard error of the sample mean shown in the parenthesis
)\small

\begin{center}

 {\begin{tabular}{cccccccc}
\hline\hline
  Method  & Model-size  & FM (\%)& $PEE\times 100$ &  $QPE\times 10^3$ 
   & $PE \times 10^3$& 
\\ \hline
 & \multicolumn{5}{c} {Correlation structure: Auto-regressive }
\\
QR           &300        & -         &1024.2(7.1)  &618.4(8.8)   &1446.4(3.1)   & \\
QR-LASSO     &39.3(1.1)  &36(0.9)  &88.7(2.5)    &60.6(2.2)    &1220.9(1.0)   &\\ 
QR-aLASSO    &17.1(0.5)  &66(1.2)  &47.9(1.9)    &37.0(2.3)     &1210.2(0.9)   &\\ 
FAL          &18.9(0.7)  &63(1.5)  &43.0(1.9)    &33.3(3.0)    &1209.5(1.0)   &\\      
FAS          &26.9(1.2)  &51(1.7)  &74.5(2.3)    &79.1(3.5)    &1223.8(1.1)   &\\
Het-QR       &8.7(0.1)   &96(0.6)  &28.8(1.3)    &24.2(1.8)    &1205.4(0.8)   &\\ 
 & \multicolumn{5}{c} {Correlation structure: Compound symmetry } 
\\
QR           &300        & -         &1041.6(7.5)  &610.9(8.3)  &1444.4(2.9)   &\\
QR-LASSO     &39.6(1.2)  &36(0.9) &86.1(2.6)    &58.0(2.4)   &1220.2(1.1)   &\\ 
QR-aLASSO    &17.0(0.5)  &66(1.3)  &46.6(2.0)    &35.7(2.3)   &1210.1(1.1)   &\\ 
FAL          &18.7(0.7)  &63(1.4)  &41.8(2.0)    &31.0(3.0)   &1208.8(1.1)   &\\      
FAS          &27.6(1.4)  &51(1.8)  &74.3(2.5)    &79.6(3.3)   &1223.5(1.1)   &\\
Het-QR       &8.7(0.1)   &96(0.7)  &27.2(1.3)    &22.3(1.9)   &1205.1(0.9)   &\\ 
 \hline
\end{tabular}}
\\
{}
\end{center}
\end{table}

\begin{table}
Table 5. Comparison of Het-QR and other methods for $p>n$  (standard error of the sample mean shown in the parenthesis
) \small  

\begin{center}

 {\begin{tabular}{cccccccc}
\hline\hline
  Method  & Model-size  & FM (\%)& $PEE \times 100$ &  $QPE \times 10^3$
   & $PE \times 10^3$ & 
\\ \hline
 & \multicolumn{5}{c} {Correlation structure: Auto-regressive } 
\\
QR-LASSO     &67.6(1.8)    &25(0.6)   &199.0(4.2)    &215.5(6.0)   &1801.4(1.7)  &\\ 
QR-aLASSO    &15.2(0.3)    &73(1.1)   &81.3(3.0)     &113.3(6.2)   &1768.8(1.6)  &\\ 
Het-QR       &10.6(0.1)    &96(0.6)  &44.6(2.5)     &61.0(7.0)    &1753.9(1.6)  &\\ 
 & \multicolumn{5}{c} {Correlation structure: Compound symmetry } 
\\
QR-LASSO     &62.9(1.8)    &27(0.7)   &182.4(4.2)   &199.8(5.9)   &1796.4(1.8)  &\\ 
QR-aLASSO    &15.2(0.3)    &74(1.0)   &76.7(2.9)    &103.5(5.9)   &1766.6(1.6)  &\\ 
Het-QR       &10.6(0.1)    &96(0.7)  &46.4(3.0)    &63.1(7.8)    &1754.2(1.7)  &\\ 
 \hline
\end{tabular}}
\\
{}
\end{center}
\end{table}

\newpage
\section{Real Data Analysis}
We collect 206 brain tumor patients each with 91 genes expression levels. All patients were de-identified. All patients were diagnosed to have glioma, one of the deadliest cancers among all cancer types. Indeed, many patients died within 1 year after the diagnosis. Glioma is associated with a number of genes. We focus on the PDGFRA gene, which encodes the alpha-type platelet-derived growth factor receptor and has been shown to be an important gene for brain tumors \citep{holland2000glioblastoma,puputti2006amplification}. We use this data set to investigate how the expression of PDGFRA is influenced by other genes.
 
For demonstration, we set $\tau$ to 0.25, 0.5, and 0.75.  For QR-LASSO, QR-aLASSO, and Het-QR, we use  cross-validation to ascertain the tuning parameter. That is, (1) we divide the data into 3 folds; (2) we use 2 folds to build the model and 1 fold to calculate the prediction errors, and this is done three times; (3) we choose the $\lambda_n$ that minimizes the prediction error as the best tuning parameter (we were not able to obtain an independent dataset with sample size $10n$ to determine the tuning parameter; it would be meaningful to compare the two procedures when such a dataset becomes available in future).  For  FAL and FAS, we follow \citep{jiang2014interquantile} to use BIC and AIC for determining the tuning parameter, and the corresponding methods are named as FAL-BIC, FAL-AIC, FAS-BIC, FAS-AIC. Hence, in total 7 approaches  are compared. We examine the  final model sizes after applying the 7 approaches. Taken all the three quantiles  together, the number of nonzero covariates of the seven models are 89 (QR-LASSO), 47 (QR-aLASSO), 25 (Het-QR), 49 (FAL-BIC), 182 (FAL-AIC), 93 (FAS-BIC), 167 (FAS-AIC). For illustration, we list some of the estimated regression coefficients in Table 6. It can be seen that the coefficients for a given gene often differ among different quantiles.  For a better view of the regression coefficients among different quantiles,  we plot the estimated coefficients for the first 30 covariates (Figure 1).  Table 6 and Figure 1 show that most models (except  FAS-BIC and FAS-AIC) demonstrate  heterogeneous sparsity, i.e., some covariates have nonzero effects in only one or two of the three quantiles. FAS-BIC and FAS-AIC do not show this type of sparsity due to the sup-norm penalty they adopt, as this penalty either selects or removes a covariate for all the quantiles. FAL-AIC and FAS-AIC models contain more nonzero estimates than FAL-BIC and FAS-BIC, consistent with the fact that BIC favors smaller models than AIC. Compared to other methods,  Het-QR yields a smaller model which may be easier to interpret and prioritize candidate genes for further functional study. 

The covariates selected by Het-QR are shown in Table 7. Consistent with the model assumptions, the estimated regression coefficients show heterogeneity among quantiles. For example, the CDKN2C gene has zero coefficient at $\tau=0.25$, and nonzero coefficients at $\tau=0.5$ and 0.75. In contrast, some other genes, such as BMP2 and SLC4A4, have nonzero coefficients across all the considered quantiles. This suggests that the expression of PDGFRA is influenced by other genes in a delicate manner that may not be fully characterized by least square methods or quantile regression methods that fail to account for the genetic heterogeneity. 
CDKN2C encodes a cyclin-dependent kinase, and BMP2 and SLC4A4 encode a bone morphogenetic protein and a sodium bicarbonate cotransporter, respectively. This indicates that PDGFRA's expression is associated with genes with a wide spectrum of cellular functions. The gene EGFR has non-positive regression coefficients, suggesting that there may be some negative control between PDGFRA and EGFR. Future biological studies may provide new insight on the gene regulation of PDGFRA. 

One main purpose of variable selection is to apply the selected variables from one dataset to other datasets to guide statistical analysis. Along this line, we further collect the brain tumor data from the cancer genome atlas (TCGA) project, which contains 567 subjects. We apply the models selected by different methods from the training data to the TCGA data to assess the prediction accuracy of the different models. 
We randomly split the TCGA data into two halves, and use one half to estimate the regression coefficients and the other half to calculate the prediction error; the prediction error is then averaged across the two halves. We repeat the random-splitting 400 times, and calculate the average of the prediction errors. Het-QR appears to have a slightly lower prediction error than the other compared ones, but the difference among the seven methods is generally small; in detail, the observed prediction errors are 1.349 (QR-LASSO), 1.351 (QR-aLASSO), 1.345 (Het-QR), 1.362 (FAL-BIC), 1.513 (FAL-AIC), 1.355 (FAS-BIC), and 1.430 (FAS-AIC).   

\newpage

\begin{figure}[!htbp]
\centering
\includegraphics[width = \textwidth,height=.8\textheight]{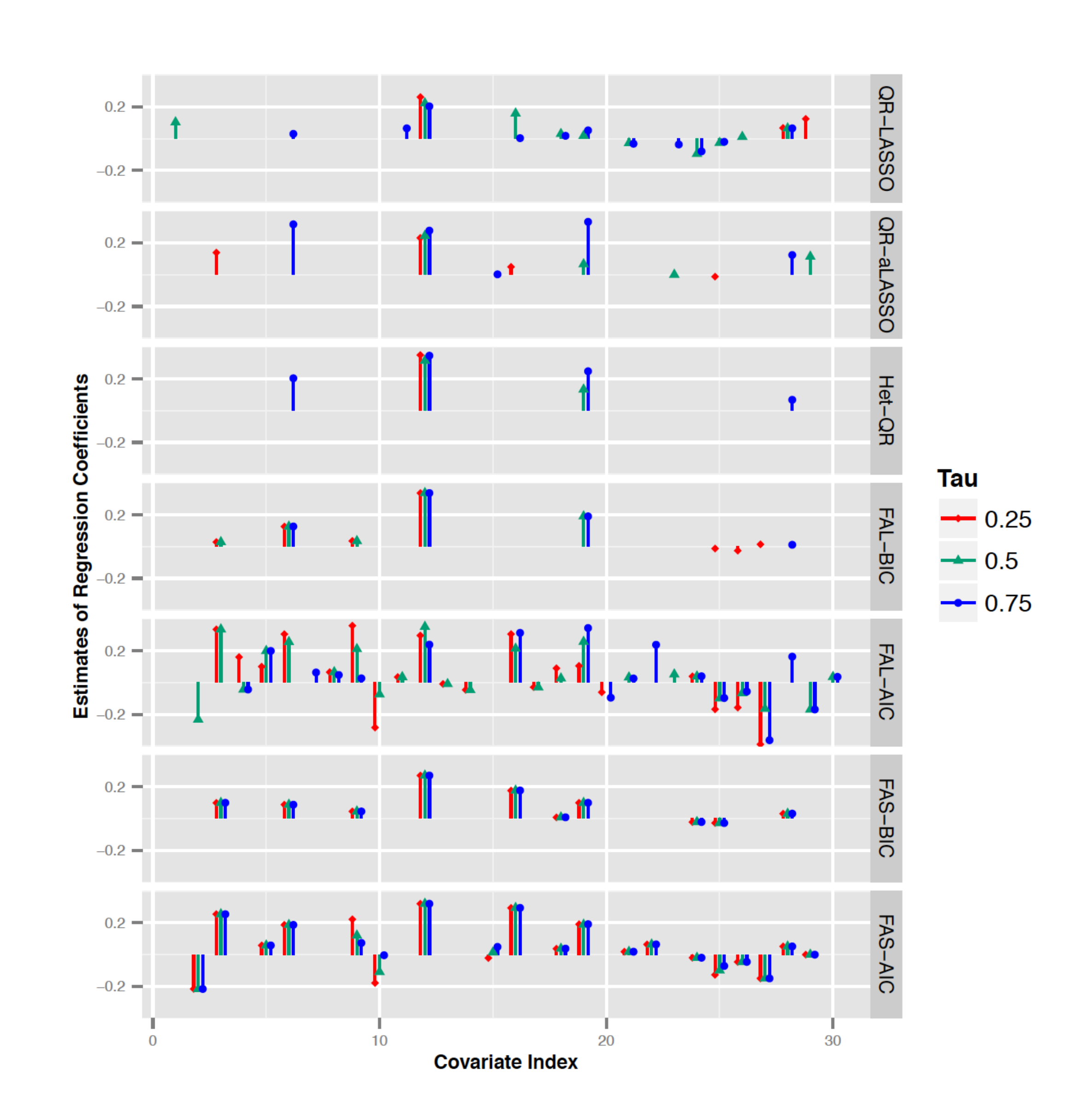} 
\label{fig1}
\caption{Graphic view of the first 30 regression coefficients estimated by different methods. 
Estimates are thresholded at 0.4 and -0.4, and only nonzero estimates are shown.    } 
\end{figure}

\begin{table}
 Table 6.  A snapshot of the estimated regression coefficients (only 5 covariates are shown) \small 
\begin{center}
{\begin{tabular}{lccccccccc} \hline ~~~Gene & $\tau$ &   QR-LASSO& QR-aLASSO
&Het-QR  &FAL-BIC &FAL-AIC &FAS-BIC &FAS-AIC &
\\
\hline 
     & 0.25&       &         &       &      &0.10     &     &0.06    &
\\
POLR2A&0.5&       &         &       &      &0.20     &     &0.06    &
\\
     & 0.75&       &         &       &      &0.20     &     &0.06    &\\
\hline
     & 0.25&        &       &       &0.13    &0.30     &0.09    &0.19    &
\\
SDHA     & 0.5&        &        &       &0.13    &0.26     &0.09    &0.19    &
\\
     & 0.75& 0.03   &0.32    & 0.2   &0.13    &0.58     &0.09    &0.19   &\\
\hline
     & 0.25&         &       &       &       &0.09     &0.01    &0.04    &
\\
CDKN2A& 0.5& 0.03    &       &       &       &0.03     &0.01    &0.04    &
\\
     & 0.75& 0.02    &       &       &       &        &0.01    &0.04    &\\
\hline
     & 0.25&         &       &       &       &0.10      &0.10    &0.19    &
\\
CDKN2C& 0.5& 0.02    &0.07   & 0.13  &0.19   &0.26     &0.10    &0.19    &
\\
     & 0.75& 0.05    &0.33   & 0.25  &0.19   &0.34     &0.10    &0.19    &\\
\hline
     & 0.25& 0.07    &       &       &       &         &0.03    &0.05    &
\\
DLL3     & 0.5& 0.06    &       &       &       &         &0.03    &0.05    &
\\
     & 0.75& 0.06    &0.12   & 0.07  &0.01   &0.16     &0.03    &0.05    &\\
\hline
\end{tabular}}
{\begin{flushleft} Note: Zero estimates are left blank.
\end{flushleft}}
\end{center}
\end{table}

\begin{table}
\begin{center}
 Table 7.   The model selected by the Het-QR method\\ \small 

 {\begin{tabular}{ccccccccc}
\hline
  &  &   \multicolumn{6}{c} {Estimated Regression Coefficients } &
\\ \cline{2-8}
Gene  &  &  & $\tau=0.25$ &  $\tau=0.5$ &  $\tau=0.75$ &
\\  
\hline
 SDHA            &  &&     &     &0.2  
\\
BMP2             &  &&0.35 &0.31 &0.34 
\\
CDKN2C           &  &&     &0.13 &0.25  
\\
DLL3             &  &&     &     &0.07 
\\
EGFR                &  &&     &-0.08 &-0.29         
\\
GRIA2                &  &&0.28 &0.22  &0.18         
\\
LTF                &  &&     &     &0.07         
\\
OLIG2                &  &&0.14 &0.30 &0.38         
\\
PLAT                &  &&0.20 &0.21 &0.25         
\\
SLC4A4                &  &&-0.21     &-0.25     &-0.24         
\\
 TAGLN                &  & & & &-0.20
\\
TMEM100                &  &&0.20 &0.17 &         
\\
 \hline
\end{tabular}}
{
}
\end{center}
\end{table}

\section{Discussion}

In this article, we have proposed a variable selection method that is able to conduct joint variable selection and estimation for multiple quantiles simultaneously. The joint selection/estimation allows one to harness the strength shared among multiple quantiles and to achieve a model that is closer to the truth. In particular, our approach is able to handle the heterogeneous sparsity, under which a covariate contributes to some (but not all) of the quantiles. By considering the heterogeneous sparsity, one can better dissect the regression structure of the trait over the covariates, which in turn leads to more accurate characterization of the underlying biological mechanism.

We have conducted a series of simulation studies to evaluate the performance of our proposed approach as well as other approaches. Our simulation studies show that the proposed method has superior performance to its peer methods. In real data analysis, our method tends to yield a sparser model than the compared methods. The benefit of achieving a sparse model is of great importance to biological studies, because it helps biological investigators to narrow down important candidate covariates (such as genes or proteins), so that research efforts can be leveraged more efficiently. 
Our analysis indicates that the regression coefficients at different quantiles can be quite heterogeneous. We suggest that the interpretation of the results be guided by biological knowledge and scientific insight, and that the variability be examined by experimental studies. FAL and FAS were mainly designed to generate interquantile shrinkage for quantile regression; when a smooth $\gamma(\tau)$ (with respect to $\tau$) is desired, these two methods are highly suitable and are indeed the only available methods to achieve such a goal. 

We have also provided theoretical proof for the proposed method under the situation that $p$ can grow to infinity.
Our exploratory experiments suggest that Het-QR can be potentially applied to `$p>n$',  although theoretical work is still needed to guide future experiments in this direction. 
Wang et al. (2012) proposed a novel approach for studying asymptotics under the `$p>n$' situation, 
and they focused on the penalties that can be written as the difference of two convex functions. The group penalty considered herein does not seem to fall into their framework.  
Further theoretical development is merited.
While we have imposed equal weights for multiple quantiles in this paper, our method can be easily extended to accommodate different weights for different quantiles. Properly chosen weights may lead to improved efficiency of the estimated parameters (Zhao and Xiao, 2013).  



\section*{Acknowledgments}

\noindent  
The authors are grateful to the AE and two reviewers for many helpful and constructive comments.
Dr. He's research is supported by the Institutional Support from the Fred Hutchinson Cancer Research Center. Dr. Kong's research is supported by Natural Sciences and Engineering Research Council of Canada. The authors thank Dr. Huixia Judy Wang for generously providing the code for FAL and FAS, which greatly expands the breadth of this manuscript.  
The authors
thank Dr. Li Hsu for helpful discussions.

The results shown here are in part based upon data generated by the TCGA Research Network: http://cancergenome.nih.gov/.

\noindent{\it Conflict of Interest}: None declared.

\section*{Appendix}

\subsection*{Transformation of the objective function}
Our proof is in vein with the proof of Proposition 1 in \cite{huang2009group}. Consider the transformed objective function
\begin{align}
&&{\min}_{  \theta,  \xi}   \sum_{m=1}^M\sum_{i=1}^{n}\,\rho_m(Y_i-U_i^{\rm T}\theta_m)+
\lambda_1\sum_{j=1}^p \xi_j +
 \sum_{j=1}^p \xi_j^{-1} \left( \sum_{m=1}^M \omega_{mj}
|\gamma_{mj}| \right) \label{AA} \\
&=&{\min}_{  \theta}  \left\{   \sum_{m=1}^M\sum_{i=1}^{n}\,\rho_m(Y_i-U_i^{\rm T}\theta_m)+ 
   \min_{\xi} \left\{
      \lambda_1\sum_{j=1}^p \xi_j +
      \sum_{j=1}^p \xi_j^{-1} \left( \sum_{m=1}^M \omega_{mj}
      |\gamma_{mj}| \right)  \right\} \right\}\notag
\end{align}

By Cauchy-Schwartz inequality, we have $$ 
      \lambda_1\sum_{j=1}^p \xi_j +
      \sum_{j=1}^p \xi_j^{-1} \left( \sum_{m=1}^M \omega_{mj}
      |\gamma_{mj}| \right)  \ge  2 \sqrt \lambda_1 \sum_{j=1}^p\left(\sum_{m=1}^M \omega_{mj}|\gamma_{mj}|\right)^{\frac{1}{2}},$$

then it follows that \eqref{AA} is equivalent to 
\begin{align}
 \min_{\theta} \sum_{m=1}^M\sum_{i=1}^{n}\,\rho_m(Y_i-U_i^{\rm T}\theta_m)
	 +  2\sqrt\lambda_1\sum_{j=1}^p\left(\sum_{m=1}^M \omega_{mj}|\gamma_{mj}|\right)^{\frac{1}{2}}. \label{BB}
\end{align}
Now, let $2\sqrt\lambda_1=n\lambda_n$, then \eqref{BB} is identical to the original objective function \eqref{0-ln}.\hfill$\Box$

%
\subsection*{Derivation of the primal and dual problem}
 Let $\lambda_{mj}=\xi_j^{-1}\omega_{mj}$, then in step 2 of Section 2.2, we aim to solve 
$${\min}_{  \gamma, \gamma_0}   \left\{ \sum_{m=1}^M\sum_{i=1}^{n}\,\rho_m(Y_i- Z_i^{\rm T}\gamma_m - \gamma_{m0}) +
 \sum_{j=1}^p  \sum_{m=1}^M \lambda_{mj}
|\gamma_{mj}| \right\}.$$ 
Let $e_n$ denote the unit vector of length $n$, and $\boldsymbol\lambda_m$ the vector of $\lambda_{mj} (j=1,\ldots, p)$. With slight abuse of notation, let $Y$ be the $n\times 1$ vector consisting of $Y_i$, and $Z$ the $n\times p$ matrix consisting of $Z_i$. The above objective function is equivalent to
$$\min_{  u_m,   v_m, \gamma_m,   s_m,   t_m} \sum_{m=1}^M \tau_m e_n^T u_m + (1-\tau_m) e_n ^T v_m + \boldsymbol\lambda_m^T s_m + \boldsymbol\lambda_m^T t_m ,$$ subject to
$u_m - v_m = Y  - Z \gamma_m - \gamma_{m0}e_n$,
$s_m - t_m= \gamma_m$, $u_m \ge 0, v_m \ge 0, s_m \ge 0,$ and $t_m \ge 0$.

Let ${\bf 0}_{r}$ be the zero vector of length $r$,  $\lambda^*=(\boldsymbol\lambda_1^T, \ldots, \boldsymbol\lambda_M^T)^T$, $s^*=(s_1^T, \ldots, s_M^T)^T$,  $t^*=(t_1^T, \ldots, t_M^T)$, $u^*=(u_1^T, \ldots, u_M^T)^T$, $v^*=(v_1^T, \ldots, v_M^T)^T$. Let $Y_{(M)}$ denote the vector in which $Y$ is stacked by $M$ times. Let  $$c =({\bf 0}_{Mp}^T, {\bf 0}_M^T, {\lambda^*}^T, {\lambda^*}^T, \tau_1 e_n^T, \ldots, \tau_M e_n^T, (1-\tau_1)e_n^T, \ldots, (1-\tau_M)e_n^T )^T,$$ $x = (\gamma^T, \gamma_{0}^T, {s^*}^T, {t^*}^T, {u^*}^T, {v^*}^T )^T,$ and $b=(Y_{(M)}^T, {\bf 0}_{Mp}^T)^T$.
Then, the linear program primal of the above objective function can be written as
$$\min_{x} c^T x,$$ subject to $Ax=b$ and $({s^*}^T, {t^*}^T, {u^*}^T, {v^*}^T )^T \ge 0$, where 
 $A$ is defined as follows. $A$ is a matrix consisting of two rows of blocks.  The first row of $A$ consists of 6 blocks, $A_{11}=I_M \otimes Z$, $A_{12}=I_M \otimes e_n$, $A_{13}=[{\bf 0}]_{nM \times Mp}$, $A_{14}=[{\bf 0}]_{nM \times Mp}$, $A_{15}=I_{nM}$, and  $A_{16}=-I_{nM}$.  The second row of $A$ consists of  the following 6 blocks, $A_{21}=I_{Mp}$, $A_{22}=[0]_{Mp\times M}$, $A_{23}=-I_{Mp}$, $A_{24}=I_{Mp}$,  $A_{25}=[{\bf 0}]_{nM \times nM}$, and $A_{26}=[{\bf 0}]_{nM \times nM}$.

 Then using standard linear program arguments, we obtain the dual as 
$$\max_{\tilde d} b^T \tilde d,$$
subject to $\tilde Z^T \tilde d = S_1 + S_2$ and $\tilde d \in [0,1]^{nM + Mp},$ where
$$S_1= \left( (1-\tau_1)e_n^T Z, \ldots, (1-\tau_M)e_n^T Z,  n(1-\tau_1),\ldots, n(1-\tau_M)  \right)^T,$$
$S_2= 1/2 \times ( R,  [{\bf 0}]_{Mp\times M} )^T e_{Mp}$, 
$R=diag( (2\boldsymbol\lambda^T, \ldots, 2\boldsymbol\lambda^T_M)^T )$, and $\tilde Z$ is defined as follows.
$\tilde Z$ consists of two rows of blocks. The first row of $\tilde Z$ includes two blocks, $\tilde Z_{11}=I_M \otimes Z$ and $\tilde Z_{12}=I_M \otimes e_n$. The second row includes two blocks, $\tilde Z_{21}=  R $ and $\tilde Z_{22}=[{\bf 0}]_{Mp\times M}$.

\subsection*{Computation time}

For $p=100$ under the auto-regressive structure, i.e., Table 2, we calculate the summary statistics of the CPU time (seconds). The average time (and the standard error of the sample mean) for QR, QR-LASSO, QR-aLASSO, FAL, FAS, Het-QR is 1.6(0.003), 10.2(0.017), 10.0(0.030), 867.0(8.059), 455.3(2.641), 90.3(0.391), respectively. For $p=600$ under the auto-regressive structure, the CPU time for QR-LASSO is 372.6(1.453); the time for QR-aLASSO and Het-QR is 362.0(1.729) and 3717.5(36.401), respectively (excluding the time for calculating penalty weights).     

\subsection*{Proof of the theorems}
We now give the proof of the theorems.
We note that throughout the proof, the upper letter $C$ in different formulas stands for different constants.


Recall that the definitions of $\mathcal{N}$, $I$, $II$, $s_n$ and $I_m$ have been given in the main text.
Define the index set
$J=\{1\leq j\leq p_n: \text{there exists}\,1\leq m\leq M\,\text{such that}\,\gamma_{mj}^* \neq 0\}$
with the cardinality $|J|=k_n$.
For every fixed $j\in J$, define the index set $M_j = \{1\leq m\leq M:\gamma_{mj}^* \neq 0\}$.
Clearly, the oracle index set $I=\{(m,j): m\in M_j, \,j\in J\}$.
For every fixed $1\leq m\leq M$, define $II_m = \{1\leq j\leq p_n:\gamma_{mj}^* = 0\}$.


We need to define some notations for our proof.
Let the vectors $\gamma_{mI}$, $\gamma_{mI}^*$ and $\hat{\gamma}_{mI}$ be the subvectors of $\gamma_m$, $\gamma_m^*$ and $\hat{\gamma}_m$ corresponding to the index set $I_m$, respectively.
Define the subvectors of $\gamma$, $\gamma^*$ and $\hat{\gamma}$ corresponding to the oracle index set $I$ as
$\gamma_I=(\gamma_{1I}^{\rm T},\cdots,\gamma_{MI}^{\rm T})^{\rm T}$,
$\gamma_I^*=(\gamma_{1I}^{*\rm T},\cdots,\gamma_{MI}^{*\rm T})^{\rm T}$ and
$\hat{\gamma}_I=(\hat{\gamma}_{1I}^{\rm T},\cdots,\hat{\gamma}_{MI}^{\rm T})^{\rm T}$.
Let $\theta_{mI} =(\gamma_{m0},\gamma_{mI}^{\rm T})^{\rm T}$, $\theta_{mI}^* =(\gamma_{m0}^*,\gamma_{mI}^{*\rm T})^{\rm T}$ and 
$\hat{\theta}_{mI} =(\hat{\gamma}_{m0},\hat{\gamma}_{mI}^{\rm T})^{\rm T}$.
Define the vector $\theta_I$ as the subvector of the parameter vector $\theta$ corresponding to the oracle index set $I$.
Recall that the vectors $\theta_I^*$ and $\hat{\theta}_I$ are the subvectors of the vectors $\theta^*$ and $\hat{\theta}$ 
corresponding to the oracle index set $I$. Clearly,
$\theta_I = (\theta_{1I}^{\rm T},\cdots,\theta_{MI}^{\rm T})^{\rm T}$,
$\theta_I^* = (\theta_{1I}^{*\rm T},\cdots,\theta_{MI}^{*\rm T})^{\rm T}$ and
$\hat{\theta}_I = (\hat{\theta}_{1I}^{\rm T},\cdots,\hat{\theta}_{MI}^{\rm T})^{\rm T}$.

Similarly, let the vectors $\gamma_{mII}$, $\gamma_{mII}^*$ and $\hat{\gamma}_{mII}$ be the subvectors of $\gamma_m$, $\gamma_m^*$ and $\hat{\gamma}_m$ corresponding to the index set $II_m$, respectively.
Define the subvectors of $\gamma$, $\gamma^*$ and $\hat{\gamma}$ corresponding to the index set $II$ as
$\gamma_{II}=(\gamma_{1II}^{\rm T},\cdots,\gamma_{MII}^{\rm T})^{\rm T}$,
$\gamma_{II}^*=(\gamma_{1II}^{*\rm T}, \cdots,\gamma_{MII}^{*\rm T})^{\rm T}$ and
$\hat{\gamma}_{II}=(\hat{\gamma}_{1II}^{\rm T}, \cdots,\hat{\gamma}_{MII}^{\rm T})^{\rm T}$.
Define the vector $\theta_{II}$ as the subvector of the parameter vector $\theta$ corresponding to the index set $II$. 
Recall that the vectors $\theta_{II}^*$ and $\hat{\theta}_{II}$ are the subvectors of the vectors $\theta^*$ and $\hat{\theta}$ corresponding to the index set $II$, respectively.
Clearly, $\theta_{II}=\gamma_{II}$, $\theta_{II}^*=\gamma_{II}^*=0$ and $\hat{\theta}_{II}=\hat{\gamma}_{II}$.

For convenience, write $\gamma=(\gamma_I^{\rm T},\gamma_{II}^{\rm T})^{\rm T}$, $\gamma^*=(\gamma_I^{*\rm T},\gamma_{II}^{*\rm T})^{\rm T}$ and  $\hat{\gamma}=(\hat{\gamma}_I^{\rm T},\hat{\gamma}_{II}^{\rm T})^{\rm T}$; and $\theta=(\theta_I^{\rm T},\theta_{II}^{\rm T})^{\rm T}$, $\theta^*=(\theta_I^{*\rm T},\theta_{II}^{*\rm T})^{\rm T}$ and $\hat{\theta}=(\hat{\theta}_I^{\rm T},\hat{\theta}_{II}^{\rm T})^{\rm T}$.

We first give a lemma related to the loss function $Q_n(\theta)$. The lemma plays an important role in the proof of our theorems.

{\sc Lemma 1.}
{\it Under conditions (L1)-(L4), we have
\begin{align*}
	Q_n(\theta)=Q_n(\theta^*)+\mathbf{A}_n^{\rm T}(\theta-\theta^*)+\frac{1}{2}(\theta-\theta^*)^{\rm T}\mathbf{B}_n(\theta-\theta^*)+R_n(\theta),
\end{align*}
where $\sup_{1\leq\,m\leq\,M}\,\sup_{\|\theta_m-\theta_m^*\|_2\,\leq\,\eta_n}\,|R_n(\theta)|= O_p(n^\frac{1}{2}p^\frac{3}{4}\eta_n^\frac{3}{2})+O_p(np^\frac{1}{2}\eta_n^3)$, and
\begin{align}\label{0-anbn}
	\mathbf{A}_n &=(A_1^{\rm T},\cdots,A_M^{\rm T})^{\rm T} \quad \text{with}\quad 
		A_m=\sum_{i=1}^n\,\bigl( I(Y_i< U_i^{\rm T}\theta_m^*) - \tau_m \bigr)\,U_i,\notag\\
	\mathbf{B}_n &=Diag(B_{11},B_{22},\cdots,B_{MM})\quad \text{with}\quad 
		B_{mm}=\sum_{i=1}^n\,f(U_i^{\rm T}\theta_m^*|\,Z_i)\,U_iU_i^{\rm T}.
\end{align}
}

{\sc Proof of Lemma 1.}\quad
Let $\psi_m(u)$ be a sub-derivative of the quantile function $\rho_m(u)$, then $\psi_m(u) =\tau_m - I(u<0)+ l_m\,I(u=0)$ with $l_m \in [-1,0]$.
Let $T_n=Q_n(\theta)-Q_n(\theta^*)$. Then, there exists an $l_{mi} \in [-1,0]$ for every $1\leq m\leq M$ and $1\leq i\leq n$ such that
\begin{align}\label{1-tn}
	T_n &= -\sum_{m=1}^M\sum_{i=1}^n\,\psi_m(Y_i-U_i^{\rm T}\,\bar{\theta}_m)\,U_i^{\rm T}(\theta_m-\theta_m^*)\notag\\
		&= -\sum_{m=1}^M\sum_{i=1}^n\,\bigl(\tau_m- I(Y_i <U_i^{\rm T}\,\bar{\theta}_m) + l_{mi}\,I(Y_i =U_i^{\rm T}\,\bar{\theta}_m)    \bigr)\,U_i^{\rm T}(\theta_m-\theta_m^*)\notag\\
		&= \sum_{m=1}^M\sum_{i=1}^n\,\bigl(I(Y_i <U_i^{\rm T}\,\theta_m^*) -\tau_m\bigr)\,U_i^{\rm T}(\theta_m-\theta_m^*)\notag\\
		&\quad -\sum_{m=1}^M\sum_{i=1}^n\,\bigl(I(Y_i <U_i^{\rm T}\,\theta_m^*) - I(Y_i <U_i^{\rm T}\,\bar{\theta}_m) \bigr)
			U_i^{\rm T}(\theta_m-\theta_m^*)\notag\\
		&\quad - \sum_{m=1}^M\sum_{i=1}^n\,l_{mi}\,I(Y_i=U_i^{\rm T}\,\bar{\theta}_m)U_i^{\rm T}(\theta_m-\theta_m^*)\notag\\
		&\equiv T_{n1} - T_{n2} - T_{n3}.	
\end{align}
where $\bar{\theta}_m$ is on the linear segment between $\theta_m$ and $\theta_m^*$, and may be written as $\bar{\theta}_m=\theta_m^*+\eta_m\,(\theta_m-\theta_m^*)$ with $\eta_m \in (0,1)$.
For $T_{n3}$, note that $Y_i$ has a continuous conditional distribution given $Z_i$, 
hence almost surely $I(Y_i=U_i^{\rm T}\,\bar{\theta}_m)=0$ for all $i=1,\cdots,n$ and $m=1,\cdots,M$, 
thus $T_{n3}=0$ almost surely.
Subsequently, we can write 
$$T_n= \bigl(T_{n1}-\mathbf{E}(T_{n1}) \bigr) - \bigl(T_{n2}-\mathbf{E}(T_{n2}) \bigr)+\mathbf{E}(T_n),$$ where $ \mathbf{E}(T_n)=\mathbf{E}(T_{n1}) +\mathbf{E}(T_{n2}) $, and $\mathbf{E}$ denotes the conditional expectation given $Z$. Note that $\mathbf{E}(T_{n1})=0$
because $\mathbf{E}(I(Y_i<U_i^{\rm T}\,\theta_m^*)-\tau_m|Z_i)=F(U_i^{\rm T}\,\theta_m^*|Z_i)-\tau_m=0$
from $\eqref{0-mod-1}$. Rename $\bigl(T_{n2}-\mathbf{E}(T_{n2})\bigr)$ as $R_{n2}$ and $\mathbf{E}(T_n)$  as $T_{n4}$, then we have
\begin{align}\label{1-tn-1}
	T_n &=  
                 T_{n1}-R_{n2}+T_{n4}.
\end{align}

For $R_{n2}$ (recall $T_{n2}$ in $\eqref{1-tn}$),
let $\zeta_{mi}(t)=I(Y_i-U_i^{\rm T}\,\theta_m^*<0) - I(Y_i-U_i^{\rm T}\,\theta_m^* < U_i^{\rm T}\,t)$, then
\begin{align}\label{1-rn2-0}
	R_{n2} 
	&= \sum_{m=1}^M \biggl( \sum_{i=1}^n\,\{ \zeta_{mi}(\eta_m(\theta_m-\theta_m^*)) -\mathbf{E}( 	
		\zeta_{mi}(\eta_m(\theta_m-\theta_m^*)))\}\,U_i \biggr)^{\rm T}(\theta_m-\theta_m^*)\notag\\
	&=\sum_{m=1}^M \phi_n(m)^{\rm T}(\theta_m-\theta_m^*)
\end{align}
where 
$\phi_n(m)= \sum_{i=1}^n\,\{\zeta_{mi}(\eta_m(\theta_m-\theta_m^*)) -\mathbf{E}( \zeta_{mi}(\eta_m(\theta_m-\theta_m^*)))\}\,U_i.$

Note that $|\zeta_{mi}(t)|=|I(Y_i-U_i^{\rm T}\,\theta_m^*<0) - I(Y_i-U_i^{\rm T}\,\theta_m^* < U_i^{\rm T}\,t)|\leq I(|Y_i-U_i^{\rm T}\,\theta_m^*|\leq |U_i^{\rm T}\,t|)$,
and that $f(t|Z_i)$ is bounded under condition (L1) and $\|U_i\|_2\leq Cp^\frac{1}{2}$ under condition (L2).
Hence,  making use of independence, under conditions (L1)and (L2), for all $1\leq m\leq M$ and $\|\theta_m-\theta_m^*\|_2\leq \eta_n$, 
we can see 
\begin{align*}  
	\mathbf{E}\|\phi_n(m)\|_2^2 
	&= \mathbf{E} \biggl\| \sum_{i=1}^n\,\{\zeta_{mi}(\eta_m(\theta_m-\theta_m^*)) -\mathbf{E}( 	
		\zeta_{mi}(\eta_m(\theta_m-\theta_m^*)))\}\,U_i \biggr\|_2^2\notag\\
	&= \sum_{i=1}^n\,\mathbf{E}\{\zeta_{mi}(\eta_m(\theta_m-\theta_m^*)) 
		-\mathbf{E}(\zeta_{mi}(\eta_m(\theta_m-\theta_m^*)))\}^2\|U_i\|_2^2 \notag\\
	&\leq Cnp\mathbf{E}\zeta_{mi}(\eta_m(\theta_m-\theta_m^*))^2
		\leq Cnp\,P\bigl(|Y_i-U_i^{\rm T}\,\theta_m^*|\leq |U_i^{\rm T}(\theta_m-\theta_m^*)|\bigr)\notag\\
        &\leq CnpP\bigl(|Y_i-U_i^{\rm T}\,\theta_m^*|\leq Cp^\frac{1}{2}\|\theta_m-\theta_m^*\|_2 )\notag\\
\end{align*}
Note that
\begin{align*}
	P\bigl(|Y_i-U_i^{\rm T}\,\theta_m^*|\leq Cp^\frac{1}{2}\|\theta_m-\theta_m^*\|_2 )
	&= F(U_i^{\rm T}\,\theta_m^*+Cp^\frac{1}{2}\|\theta_m-\theta_m^*\|_2 )
		-F(U_i^{\rm T}\,\theta_m^*-Cp^\frac{1}{2}\|\theta_m-\theta_m^*\|_2  ) \\
	&\leq \,|f(\xi_{mi}|Z_i)|p^\frac{1}{2}\|\theta_m-\theta_m^*\|_2	
		\leq C p^\frac{1}{2}\|\theta_m-\theta_m^*\|_2, 
\end{align*}
where $\xi_{mi}$ is between $U_i^{\rm T}\,\theta_m^*+Cp^\frac{1}{2}\|\theta_m-\theta_m^*\|_2$ and $U_i^{\rm T}\,\theta_m^*-Cp^\frac{1}{2}\|\theta_m-\theta_m^*\|_2$.
Thus, 
$
\mathbf{E}\|\phi_n(m)\|_2^2 \le Cnp^\frac{3}{2}\eta_n.
$
By Chebyshev's inequality, we get
$\sup_{1\leq\,m\leq\,M}\sup_{\|\theta_m-\theta_m^*\|_2\,\leq\,\eta_n}\,\|\phi_n(m)\|_2 =O_p(n^\frac{1}{2}p^\frac{3}{4}\eta_n^\frac{1}{2})$.
Together with $\eqref{1-rn2-0}$, by Cauchy-Schwartz's inequality, we get 
\begin{align}\label{1-rn2}
	\sup_{1\leq\,m\leq\,M}\,\sup_{\|\theta_m-\theta_m^*\|_2\,\leq\,\eta_n}\,|R_{n2}| 
	= O_p(n^{\frac{1}{2}}p^{\frac{3}{4}}{\eta_n}^{\frac{1}{2}}) O(\eta_n )
	= O_p(n^{\frac{1}{2}}p^{\frac{3}{4}}{\eta_n}^{\frac{3}{2}}).
\end{align}
For $T_{n4}$, the third term  in $\eqref{1-tn-1}$,  write $\mathbf{E}(T_n)=e_n(\theta)-e_n(\theta^*)$, where
\begin{align*}
	e_n(\theta)\equiv  \sum_{m=1}^M\sum_{i=1}^n\,\mathbf{E}\,\rho_m\bigl(Y_i-U_i^{\rm T}\theta_m\bigr).
\end{align*}
 $\mathbf{E}\,\rho_m\bigl(Y_i-U_i^{\rm T}\theta_m\bigr)$ is second order differentiable with respect to $\theta_m$ under condition (L1), with gradient  
$G_{mi}(\theta)= -\mathbf{E}\,\bigl\{ (\tau_m-I(Y_i<U_i^{\rm T}\theta_m))U_i \bigr\}=\bigl(F(U_i^{\rm T}\theta_m\bigr|Z_i)-\tau_m\bigr)\,U_i$,
and Hessian matrix $H_{mi}(\theta)= f(\bigl(U_i^{\rm T}\theta_m)\bigr|Z_i\bigr)\,U_iU_i^{\rm T}$. 

Let $G(\theta)$ and $H(\theta)$ be gradient and Hessian matrix of $e_n(\theta)$, then 
$G(\theta)= \sum_{m=1}^M\sum_{i=1}^n\,G_{mi}(\theta)$ and 
$H(\theta)= \sum_{m=1}^M\sum_{i=1}^n\,H_{mi}(\theta) = \sum_{m=1}^M\sum_{i=1}^n\,f(\bigl(U_i^{\rm T}\theta_m)\bigr|Z_i\bigr)\,U_iU_i^{\rm T}$.
It is easy to see that $G(\theta^*)=0$ by $F(U_i^{\rm T}\theta_m^*|Z_i)=\tau_m$ in $\eqref{0-mod-1}$.
By Taylor expansion of $T_{n4}=\mathbf{E}(T_{n})= e_n(\theta)-e_n(\theta^*)$ at $\theta^*$, we have
\begin{align}\label{1-tn4}
	T_{n4}
	&=\frac{1}{2}(\theta-\theta^*)^{\rm T}\,H\bigl(\theta^*+\xi(\theta-\theta^*)\bigr)(\theta-\theta^*)\notag\\
	&=\frac{1}{2}\sum_{m=1}^M\sum_{i=1}^n\,f(\zeta_{mi}|Z_i)\,
			(\theta_m-\theta_m^*)^{\rm T}\,U_iU_i^{\rm T}(\theta_m-\theta_m^*),\notag\\
	&=\frac{1}{2}\sum_{m=1}^M\sum_{i=1}^n\,\bigl\{f(\zeta_{mi}|Z_i)-f(U_i^{\rm T}\theta_m^*|Z_i)\bigr\}\,
			(\theta_m-\theta_m^*)^{\rm T}\,U_iU_i^{\rm T}(\theta_m-\theta_m^*)\notag\\
	&\quad+\frac{1}{2}\sum_{m=1}^M\sum_{i=1}^n\,f(U_i^{\rm T}\theta_m^*|Z_i)(\theta_m-\theta_m^*)^{\rm T}U_iU_i^{\rm T}
            (\theta_m-\theta_m^*)\notag\\
	&\equiv R_{n4}+T_{n41}.			
\end{align}
where $\xi\in(0,1)$ and $\zeta_{mi}$ is between $U_i^{\rm T}\theta_m$ and $U_i^{\rm T}\theta_m^*$.
Trivially,
\begin{align}\label{1-rn4}
	|R_{n4}| \leq C\sum_{m=1}^M\sup_{1\leq i\leq n}\bigl|f(\zeta_{mi}|Z_i)-f(U_i^{\rm T}\theta_m^*|Z_i)\bigr|
	(\theta_m-\theta_m^*)^{\rm T}\sum_{i=1}^n\,\,U_iU_i^{\rm T}(\theta_m-\theta_m^*).
\end{align}
Note that $f'(t|Z_i)$ is bounded under condition (L1), $\|U_i\|_2\leq Cp^\frac{1}{2}$ under condition (L2), 
and $\lambda_{max}\bigl(n^{-1}\sum_{i=1}^n\,U_iU_i^{\rm T}\bigr)\leq C$ under condition (L3).
Hence, for all $1\leq i\leq n$, $1\leq m\leq M$, and $\|\theta_m-\theta_m^*\|_2\leq \eta_n$, we have
$|f(\zeta_{mi}|Z_i)-f(U_i^{\rm T}\theta_m^*|Z_i)|\leq C|U_i^{\rm T}(\theta_m-\theta_m^*)|\leq\,C\|U_i\|_2\|\theta_m-\theta_m^*\|_2\leq\,Cp^\frac{1}{2}\eta_n$,
and $(\theta_m-\theta_m^*)^{\rm T}\sum_{i=1}^n\,\,U_iU_i^{\rm T}(\theta_m-\theta_m^*)\leq n\lambda_{max}\bigl(n^{-1}\sum_{i=1}^n\,U_iU_i^{\rm T}\bigr)\|\theta_m-\theta_m^*\|_2^2\leq Cn\eta_n^2$.
Hence, from $\eqref{1-rn4}$, we get 
$\sup_{1\leq\,m\leq\,M}\,\sup_{\|\theta_m-\theta_m^*\|_2\,\leq\,\eta_n}\,|R_{n4}| \leq Cnp^\frac{1}{2}\eta_n^3$.
Together with $\eqref{1-tn-1}$, $\eqref{1-rn2}$ and $\eqref{1-tn4}$, we obtain
$T_n=T_{n1}+T_{n41}+R_n(\theta)$,
where $\sup_{1\leq\,m\leq\,M}\,\sup_{\|\theta_m-\theta_m^*\|_2\,\leq\,\eta_n}\,|R_{n}(\theta)|
	= O_p(n^{\frac{1}{2}}p^{\frac{3}{4}}{\eta_n}^{\frac{3}{2}})+O_p(np^\frac{1}{2}\eta_n^3)$, and 
\begin{align*}
	T_{n1} 	&= \sum_{m=1}^M\sum_{i=1}^n\,\bigl(I(Y_i <U_i^{\rm T}\theta_m^*) -\tau_m\bigr)\,U_i^{\rm T}(\theta_m-\theta_m^*)
			= \mathbf{A}_n^{\rm T}(\theta-\theta^*),\notag\\
	T_{n41} &= \frac{1}{2}\sum_{m=1}^M\sum_{i=1}^n\,f(U_i^{\rm T}\theta_m^*|Z_i)
				(\theta_m-\theta_m^*)^{\rm T}U_iU_i^{\rm T}(\theta_m-\theta_m^*)
			= \frac{1}{2}(\theta-\theta^*)^{\rm T}\mathbf{B}_n(\theta-\theta^*).
\end{align*}
 This completes the proof of the lemma.\hfill$\Box$


{\sc Proof of Theorem 1.}\quad
Recall the definition of $L_n(\theta)= Q_n(\theta)+P_n(\gamma)$ in $\eqref{0-ln}$.
Let $\theta-\theta^*=\nu_n\,u$, where $\nu_n>0$, $u\in \mathbb{R}^{M(p+1)}$ and $\|u\|_2=1$. It is easy to see that
$\|\theta-\theta^*\|_2=\nu_n$.
Based on the continuity of $L_n$, if we can prove that in probability
\begin{align}\label{2-1}
	\inf_{\|u\|_2=1}\,L_n(\theta^*+\nu_nu)>L_n(\theta^*),
\end{align}
then the minimal value point of $L_n(\theta^*+\nu_nu)$ on $\{u:\|u\|_2\leq1\}$ exists and lies in the unit ball $\{u:\|u\|_2\leq1\}$ in probability. We will prove that $\eqref{2-1}$ holds.

For $Q_n(\theta)$, because $\|\theta_m-\theta_m^*\|_2\,\leq\,\|\theta-\theta^*\|_2$ for all $m=1,2,\cdots,M$,
by Lemma 1 with $\eta_n=\nu_n$ under conditions (L1) to (L4), we get
\begin{align}\label{2-qn}	
	q_n(\theta) \equiv Q_n(\theta)-Q_n(\theta^*)
	=\mathbf{A}_n^{\rm T}(\theta-\theta^*)+\frac{1}{2}(\theta-\theta^*)^{\rm T}\mathbf{B}_n(\theta-\theta^*)+R_n(\theta),
\end{align}
where 
$\sup_{\|\theta-\theta^*\|_2\leq \nu_n}\,|R_n(\theta)|= O_p(n^{\frac{1}{2}}p^{\frac{3}{4}}{\nu_n}^{\frac{3}{2}})+O_p(np^\frac{1}{2}\nu_n^3)$.

For $P_n(\gamma)$, let $p_n(\gamma) = P_n(\gamma)-P_n(\gamma^*)$.
Define $p_{1n}(\gamma)= p_n(\gamma_I,0)$, that is, 
\begin{align}\label{2-p1n}
	p_{1n}(\gamma) =  n\lambda_n\,\sum_{j\in J}\,\bigl(\sum_{m\in M_j}\,\omega_{mj}|\gamma_{mj}|\bigr)^{\frac{1}{2}}
		- n\lambda_n\,\sum_{j\in J}\,\bigl(\sum_{m\in M_j}\,\omega_{mj}|\gamma_{mj}^*|\bigr)^{\frac{1}{2}}.
\end{align}
Clearly, $p_{1n}(\gamma) \leq p_n(\gamma)$ and $p_{1n}(\gamma^*)= p_n(\gamma^*)=0$.

Define $l_n(\theta)=q_n(\theta)+p_n(\gamma)$ and $l_{1n}(\theta)=q_n(\theta)+p_{1n}(\gamma)$, both of which are continuous.
Clearly, $l_{1n}(\theta)\leq l_n(\theta)$ and $l_{1n}(\theta^*)=l_n(\theta^*)=0$.
Note that $\eqref{2-1}$ is equivalent to that in probability
\begin{align}\label{2-2}
	\inf_{\|u\|_2=1}\,l_{1n}(\theta^*+\nu_nu)>0.
\end{align}
Note that
\begin{align*}
	|p_{1n}(\gamma)| 
	&= n\lambda_n\,\bigl|\sum_{j\in J}\,\bigl(\sum_{m\in M_j}\,\omega_{mj}|\gamma_{mj}|\bigr)^{\frac{1}{2}}
	 -\sum_{j\in J}\,\bigl(\sum_{m\in M_j}\,\omega_{mj}|\gamma_{mj}^*|\bigr)^{\frac{1}{2}} \bigr|\notag\\
	&\leq n\lambda_n\,\sum_{j\in J}\,\bigl|\bigl(\sum_{m\in M_j}\,\omega_{mj}|\gamma_{mj}|\bigr)^{\frac{1}{2}}
		 -\bigl(\sum_{m\in M_j}\,\omega_{mj}|\gamma_{mj}^*|\bigr)^{\frac{1}{2}} \bigr|
	\leq n\lambda_n\,\sum_{j\in J}\,\bigl(\sum_{m\in M_j}\,\omega_{mj}|\gamma_{mj} - \gamma_{mj}^* |\bigr)^{\frac{1}{2}},
\end{align*}
where the last inequality follows from the fact that $|\sqrt{|x|}-\sqrt{|y|}|\leq \sqrt{|x-y|}$.
Note that $\gamma_{mj}-\gamma_{mj}^* = \nu_n\,u_{mj}$ where $u_{mj}$ is a component of $u$. Hence,
\begin{align}\label{2-tpn}
	|p_{1n}(\gamma)|
	&\leq n\lambda_n\,d_{nI}^\frac{1}{2}\nu_n^\frac{1}{2}\,\sum_{j\in J}\,\bigl(\sum_{m\in M_j}\,|u_{mj}|\bigr)^{\frac{1}{2}}
	\leq n\lambda_n\,d_{nI}^\frac{1}{2}\nu_n^\frac{1}{2}k_n^\frac{1}{2}\bigl(\sum_{(m,j)\in I}\,|u_{mj}|\bigr)^{\frac{1}{2}}\notag\\
	&\leq n\lambda_nd_{nI}^\frac{1}{2}\nu_n^\frac{1}{2}k_n^\frac{1}{2}s_n^\frac{1}{4}\|u\|_2^{\frac{1}{2}}
	=n\lambda_nd_{nI}^\frac{1}{2}\nu_n^\frac{1}{2}s_n^\frac{3}{4}\|u\|_2^{\frac{1}{2}},
\end{align}
where $k_n=|J|\leq s_n$.
By $\eqref{2-qn}$ and the above, we have
\begin{align}\label{2-ln-1}
	l_{1n}(\theta^*+\nu_nu) = q_n(\theta^*+\nu_n u)+p_{1n}(\theta^*+\nu_n u)
	= \nu_n\mathbf{A}_n^{\rm T}u+\frac{1}{2}\nu_n^2 u^{\rm T} \mathbf{B}_n u + R_L(u),
\end{align}
where $\sup_{\|u\|_2\,=\,1}\,|R_L(u)|= O_p(n^{\frac{1}{2}}p^{\frac{3}{4}}{\nu_n}^{\frac{3}{2}})+O_p(np^\frac{1}{2}\nu_n^3)
	 +O(n\lambda_n\,d_{nI}^\frac{1}{2}\nu_n^\frac{1}{2}s_n^\frac{3}{4})$.

For the quadratic term in $\eqref{2-ln-1}$, from condition (L5), 
\begin{align}\label{2-quad-1}
	\frac{1}{2}\nu_n^2 u^{\rm T} \mathbf{B}_n u \geq \frac{C_4}{2}n\nu_n^2\|u\|_2^2.
\end{align}
 For the linear term $\mathbf{A}_n^{\rm T}u$ in $\eqref{2-ln-1}$, by the independence of $(Z_i,Y_i)$ and $(Z_j,Y_j)$ for all $i\neq j$ 
and the fact that $E(I(Y_i< U_i^{\rm T}\theta_m^*) - \tau_m|Z_i)=0$, we get
\begin{align*}
	E(A_n^{\rm T}A_n)
	&= E\sum_{m=1}^M\,\sum_{i=1}^n\,\bigl( I(Y_i< U_i^{\rm T}\theta_m^*) - \tau_m \bigr)\,U_i^{\rm T} 
		\sum_{j=1}^n\,\bigl( I(Y_j< U_j^{\rm T}\theta_m^*) - \tau_m \bigr)\,U_j\notag\\
	&=\sum_{m=1}^M\,\sum_{i=1}^n\,E\bigl(\bigl( I(Y_i< U_i^{\rm T}\theta_m^*) - \tau_m \bigr)^2\| U_i\|_2^2\bigl)
	\leq Cnp.
\end{align*}
Then, it follows that 
\begin{align}\label{2-an-2}
	\|\mathbf{A}_n\|_2=O_p((np)^{\frac{1}{2}}),
\end{align} 
which implies that $\sup_{\|u\|_2\leq 1}|\mathbf{A}_n^{\rm T}u|= O_p((np)^{\frac{1}{2}})$.
Together with $\eqref{2-ln-1}$ and $\eqref{2-quad-1}$, in probability,
\begin{align}\label{2-order-2-2}
	\inf_{\|u\|_2=1}\,l_{1n}(\theta^*+\nu_nu)
	&\geq \frac{C_4}{2}n\nu_n^2
			-C(np)^{\frac{1}{2}}\nu_n - Cn^{\frac{1}{2}}p^{\frac{3}{4}}{\nu_n}^{\frac{3}{2}}-Cnp^{\frac{1}{2}}\nu_n^3 
			-Cn\lambda_n\,d_{nI}^\frac{1}{2}\nu_n^\frac{1}{2}s_n^\frac{3}{4}\notag\\
	&\geq \frac{C_4}{2}n\nu_n \bigl\{\nu_n
			-Cn^{-\frac{1}{2}}p^{\frac{1}{2}} - Cn^{-\frac{1}{2}}p^{\frac{3}{4}}{\nu_n}^{\frac{1}{2}}-Cp^{\frac{1}{2}}\nu_n^2 
			- C\lambda_nd_{nI}^\frac{1}{2}s_n^\frac{3}{4}\nu_n^{-\frac{1}{2}} \bigr\}.
\end{align}
Now take $\nu_n=C_0(n^{-{\frac{1}{2}}}p^{\frac{1}{2}})$ where $C_0$ is a sufficiently large constant.
Under condition (L4), i.e., $\lambda_nd_{nI}^\frac{1}{2}=o(s_n^{-\frac{3}{4}}p_n^{\frac{3}{4}}n^{-\frac{3}{4}})$,
for the last three terms in $\eqref{2-order-2-2}$, we can check that
	$n^{-\frac{1}{2}}p^{\frac{3}{4}}{\nu_n}^{\frac{1}{2}} \leq Cn^{-\frac{1}{4}+\frac{1}{2}\alpha_1}\nu_n=o(\nu_n)$,
	$p^\frac{1}{2}\nu_n^2 \leq Cn^{-\frac{1}{2}+\alpha_1}\nu_n=o(\nu_n)$,
	and $\lambda_n\,d_{nI}^\frac{1}{2}s_n^\frac{3}{4}\nu_n^{-\frac{1}{2}} \leq C\lambda_n\,d_{nI}^\frac{1}{2}s_n^{\frac{3}{4}}p_n^{-\frac{3}{4}}n^\frac{3}{4}\nu_n=o(\nu_n)$.
Hence, 
\begin{align}\label{2-order-3}
	\inf_{\|u\|_2=1}\,l_{1n}(\theta^*+\nu_nu)\geq Cn\nu_n^2\to \infty\quad\text{in probability}.
\end{align}
Therefore, in probability there exists a local minimizer $\hat{\theta}$ of $L_n(\theta)$ such that $\|\hat{\theta}-\theta^*\|_2<\nu_n$.
This completes the proof of the theorem.\hfill$\Box$


{\sc Proof of Theorem 2.}\quad
For the quantile function $Q_n(\theta)$, because $\|\theta_m-\theta_m^*\|_2\,\leq\,\|\theta-\theta^*\|_2$ for all $1\leq\,m\leq\,M$,
by Lemma 1 with $\eta_n=\nu_n$, under conditions (L1)-(L4), we have
\begin{align}\label{3-qn-0}
	Q_n(\theta)-Q_n(\theta^*) = \mathbf{A}_n^{\rm T}(\theta-\theta^*)+\frac{1}{2}(\theta-\theta^*)^{\rm T}\mathbf{B}_n(\theta-\theta^*)+R_n(\theta),
\end{align}
where $\sup_{\|\theta-\theta^*\|_2\,\leq\,\nu_n}\,|R_n(\theta)|
	=O_p(n^{\frac{1}{2}}p^{\frac{3}{4}}{\nu_n}^{\frac{3}{2}})+O_p(np^\frac{1}{2}\nu_n^3)$.

Let $\theta-\theta^* = \nu_n u$ where $\nu_n>0$ and $u\in \mathbb{R}^{M(p+1)}$. 
Then $\|\theta-\theta^*\|_2\leq \nu_n$ if and only if $\|u\|_2\leq 1$. 
Let $u=(u_I^{\rm T},u_{II}^{\rm T})^{\rm T}$, where $u_I$ and $u_{II}$ are are subvectors of $u$ corresponding to the index sets $I$ and $II$, respectively.
Clearly, $\|u\|_2^2=\|u_I\|_2^2+\|u_{II}\|_2^2$. 
Note that $\theta_I=\theta_I^* +\nu_n u_I$ and $\theta_{II}=\theta_{II}^* +\nu_n u_{II}=\nu_n u_{II}$.

Define the ball $\tilde{\Theta}_n=\{\theta=\theta^*+\nu_n u \in \Theta_n :\,\|u\|_2\leq 1\}$ with $\nu_n=C_0n^{-\frac{1}{2}}p^{\frac{1}{2}}$.
For any $\theta = (\theta_I^{\rm T},\theta_{II}^{\rm T})^{\rm T} \in  \tilde{\Theta}_n$, we can see  
$\|(\theta_I^{\rm T},0^{\rm T})^{\rm T}-\theta^*\|_2= \|\theta_I-\theta_I^*\|_2\leq \nu_n$, and
$\|\theta_{II}\|_2 = \nu_n\|u_{II}\|_2$, where $\|u_{II}\|_2\leq 1$.

Consider that 
$Q_n(\theta_I,\theta_{II}) - Q_n(\theta_I,0)= Q_n(\theta_I,\theta_{II}) - Q_n(\theta^*) -\bigl( Q_n(\theta_I,0) - Q_n(\theta^*)\bigr)$
\begin{align*}
	&=\frac{1}{2}(0^{\rm T},\theta_{II}^{\rm T})\mathbf{B}_n(0^{\rm T},\theta_{II}^{\rm T})^{\rm T} + \mathbf{A}_{n}^{\rm T}(0^{\rm T},\theta_{II}^{\rm T})^{\rm T}
		-(\theta_I^{\rm T}-\theta_I^{*\rm T},0^{\rm T})\mathbf{B}_n(0^{\rm T},\theta_{II}^{\rm T})^{\rm T} + \bigl( R_n(\theta)-R_n(\theta_I,0) \bigr)\notag\\
	&\equiv I_{n1}+I_{n2}+I_{n3}+r_{1n}(\theta)
\end{align*}
From $\eqref{2-an-2}$, we can see 
$\sup_{\theta \in \tilde{\Theta}_n}\,|I_{n2}| \leq \|\mathbf{A}_n\|_2\|\theta_{II}\|_2 = O_p((np)^{\frac{1}{2}})\nu_n\|u_{II}\|_2= O_p(p\|u_{II}\|_2)$.
Under condition (L5), 
$\|\mathbf{B}_n(0^{\rm T},\theta_{II}^{\rm T})^{\rm T}\|_2^2= (0^{\rm T},\theta_{II}^{\rm T})\mathbf{B}_n^2(0^{\rm T},\theta_{II}^{\rm T})^{\rm T}
\leq n^2\lambda_{\max}^2(n^{-1}B_n)\|\theta_{II}\|_2^2\leq Cn^2\nu_n^2\|u_{II}\|_2^2$, 
we have
$\sup_{\theta\in \tilde{\Theta}_n}\,|I_{n3}|\leq \|\theta_I-\theta_I^*\|_2\|\mathbf{B}_n(0^{\rm T},\theta_{II}^{\rm T})^{\rm T} \|_2 \leq Cn\nu_n^2\leq Cp$.
From $\eqref{3-qn-0}$, we have 
$\sup_{\theta\in \tilde{\Theta}_n}\,|r_{1n}(\theta)|=O_p(n^{-\frac{1}{4}}p^{\frac{3}{2}})$.
Hence, 
\begin{align}\label{3-qn}
	Q_n(\theta_I,\theta_{II}) - Q_n(\theta_I,0)= I_{n1} + r_{2n}(\theta),
\end{align}
where $\sup_{\theta\in \tilde{\Theta}_n}\,|r_{2n}(\theta)|=O_p(p\|u_{II}\|_2)+O_p(n^{-\frac{1}{4}}p^{\frac{3}{2}})$.
Recall $L_n(\theta_I,\theta_{II}) - L_n(\theta_I,0)=Q_n(\theta_I,\theta_{II}) - Q_n(\theta_I,0) + P_n(\gamma_I,\gamma_{II})-P_n(\gamma_I,0)$.
From $\eqref{3-qn}$, we get 
\begin{align}\label{3-ln-0}
	L_n(\theta_I,\theta_{II}) - L_n(\theta_I,0) \geq  P_n(\gamma_I,\gamma_{II})-P_n(\gamma_I,0)+ r_{2n}(\theta).
\end{align}
Note that $(n\lambda_n)^{-1}\bigl(P(\gamma_I,\gamma_{II})-P_n(\gamma_I,0)\bigr)$
\begin{align}\label{3-pn-1}
	&= \sum_{j=1}^p\,\bigl(\sum_{m=1}^M\,\omega_{mj}|\gamma_{mj}|\bigr)^{\frac{1}{2}} 
		-\sum_{j\in J}\,\bigl(\sum_{m\in M_j}\,\omega_{mj}|\gamma_{mj}|\bigr)^{\frac{1}{2}}\notag\\
	&= \sum_{j\in J^c}\,\bigl(\sum_{m=1}^M\,\omega_{mj}|\gamma_{mj}|\bigr)^{\frac{1}{2}} 
		+\sum_{j\in J}\,\bigl(\sum_{m=1}^M\,\omega_{mj}|\gamma_{mj}|\bigr)^{\frac{1}{2}} 
		-\sum_{j\in J}\,\bigl(\sum_{m\in M_j}\,\omega_{mj}|\gamma_{mj}|\bigr)^{\frac{1}{2}}\notag\\
	&\geq \sum_{j\in J^c}\frac{\sum_{m=1}^M\,\omega_{mj}|\gamma_{mj}|}{2\bigl(\sum_{m=1}^M\,\omega_{mj}|\gamma_{mj}|\bigr)^{\frac{1}{2}}}
		+\sum_{j\in J}\frac{\sum_{m\in M_j^c}\omega_{mj}|\gamma_{mj}|}{2\bigl(\sum_{m=1}^M\omega_{mj}|\gamma_{mj}|\bigr)^{\frac{1}{2}}}.
\end{align}
For all $\theta \in \tilde{\Theta}_n$, we have $|\gamma_{mj}|\leq |\gamma_{mj}^*|+1\leq C$, which implies that
$\sum_{m=1}^M\omega_{mj}|\gamma_{mj}|\leq C\|\omega_n\|_\infty$ for all $j=1,2,\cdots,p$, 
where $\|\omega_n\|_\infty =\max_{1\leq m\leq M,\,1\leq j\leq p}|\omega_{mj}|$.
Recall that $d_{nII}=\min_{(m,j)\in II}\{\omega_{mj}\}\|\omega_n\|_\infty ^{-\frac{1}{2}}$.
From $\eqref{3-pn-1}$, it follows that for all $\theta\in \tilde{\Theta}_n$,
\begin{align}\label{3-pn-0}
	P_n(\gamma_I,\gamma_{II})-P_n(\gamma_I,0)
	&\geq Cn\lambda_n\,\|\omega_n\|_\infty ^{-\frac{1}{2}}
		\bigl(\sum_{j\in J^c}\sum_{m=1}^M\,\omega_{mj}|\gamma_{mj}| + \sum_{j\in J}\sum_{m\in M_j^c}\omega_{mj}|\gamma_{mj}|\bigr)\notag\\
	&= Cn\lambda_n\,\|\omega_n\|_\infty ^{-\frac{1}{2}}\sum_{(m,j)\in II}\,\omega_{mj}|\gamma_{mj}|
		\geq Cn\lambda_n\,d_{nII} \sum_{(m,j)\in II}|\gamma_{mj}|\notag\\
	&\geq Cn\lambda_n\,d_{nII} \|\gamma_{II}\|_2 = C\lambda_n\,d_{nII}(np)^{\frac{1}{2}}\|u_{II}\|_2.	
\end{align}

Define $\Omega_n = \{\theta=\theta^*+\nu_n u \in \tilde{\Theta}_n :\|u_{II}\|_2 > 0 \}$
and $\Omega_n^c = \{\theta=\theta^*+\nu_n u \in \tilde{\Theta}_n :u_{II} = 0\}$.
Clearly, $\tilde{\Theta}_n=\Omega_n \cup \Omega_n^c $.
From $\eqref{3-ln-0}$ and $\eqref{3-pn-0}$, we obtain in probability 
\begin{align*}
	\inf_{\theta\in \Omega_n}\,\bigl(L_n(\theta_I,\theta_{II}) - L_n(\theta_I,0)\bigr)
	&\geq \inf_{\theta\in \Omega_n}\,\bigl(P_n(\gamma_I,\gamma_{II})-P_n(\gamma_I,0)\bigr)
			-\sup_{\theta\in \tilde{\Theta}_n}\,|r_{2n}(\theta)|\notag\\
	&\geq \tilde C_1\lambda_n\,d_{nII}(np)^{\frac{1}{2}}\|u_{II}\|_2 - \tilde C_2 p\|u_{II}\|_2 - \tilde C_2 n^{-\frac{1}{4}}p^{\frac{3}{2}}\\
	&\geq p\biggl( \|u_{II}\|_2\bigl(\tilde C_1\lambda_n\,d_{nII}n^{\frac{1}{2}}p^{-\frac{1}{2}} -\tilde C_2\bigr) - \tilde C_2 n^{-\frac{1}{4}}p^{\frac{1}{2}}\biggr).
\end{align*}
where $\tilde C_1$ and $\tilde C_2$ are positive constants.
Under the given conditions, $\lambda_n\,d_{nII}n^{\frac{1}{2}}p^{-\frac{1}{2}}\to\infty $ and $n^{-\frac{1}{4}}p^{\frac{1}{2}}\to 0$
as $n\to \infty$. 
Hence,  
$\inf_{\theta\in \Omega_n}\,\bigl(L_n(\theta_I,\theta_{II}) - L_n(\theta_I,0)\bigr)>0$ in probability.
Thus,
$\inf_{\theta\in \Omega_n}\,L_n(\theta) \geq \inf_{\theta\in \Omega_n}\,\bigl(L_n(\theta_I,\theta_{II}) - L_n(\theta_I,0)\bigr)
+\inf_{\theta\in \Omega_n}L_n(\theta_I,0) > \inf_{\theta\in \Omega_n}L_n(\theta_I,0) 
= \inf_{\theta\in \tilde{\Theta}_n} L_n(\theta_I,0) \geq \inf_{\theta\in \tilde{\Theta}_n} L_n(\theta)$,
which implies that
$\inf_{\theta\in\tilde{\Theta}_n}\,L_n(\theta)=\inf_{\theta\in \Omega_n^c}\,L_n(\theta)$. Therefore,
the minimal value point of $L_n(\theta)$ on $\tilde{\Theta}_n$ only lies in its subset $\Omega_n^c$.

From Theorem 1, we know that in probability $\hat{\theta} \in \tilde{\Theta}_n$ and that $\hat{\theta}$ is a local minimizer of $L_n(\theta)$. 
Hence, $\hat{\theta} \in \Omega_n^c$ in probability, which implies that $\hat{\gamma}_{II} = 0$.
\hfill$\Box$


{\sc Proof of Theorem 3.}\quad
Let $\theta_I-\theta_I^*=\nu_n u$ where $\nu_n>0$ and $u\in \mathbb{R}^{M+s_n}$.
Because $\|\theta_{mI}-\theta_{mI}^*\|_2\,\leq\,\|\theta_I-\theta_I^*\|_2 = \nu_n\|u\|_2$ for all $1\leq\,m\leq\,M$,
and due to conditions (L1)-(L3) and that $0<\alpha_0<\alpha_1<\frac{1}{6}$, Lemma 1 implies that
\begin{align}\label{4-qn}
	Q_n(\theta_I,0)-Q_n(\theta_I^*,0) =	Q_n(\theta_I^*+\nu_n u,0)-Q_n(\theta_I^*,0)
	=\nu_n\mathbf{A}_{nI}^{\rm T}u + \frac{1}{2}\nu_n^2u^{\rm T}\mathbf{B}_{nI}u+r_n(u),
\end{align}
where $\sup_{\|u\|_2\,\leq\,1 }\,|r_n(u)|= O_p(n^{\frac{1}{2}}s_n^{\frac{3}{4}}\nu_n^{\frac{3}{2}})+O_p(ns_n^{\frac{1}{2}}\nu_n^3)$,
$\mathbf B_{nI}$ is given in $\eqref{4-bn}$,
\begin{align}\label{4-an}
	\mathbf{A}_{nI} &=(A_1^{\rm T}, \cdots,A_M^{\rm T})^{\rm T} \quad \text{with}\quad 
		A_m=\sum_{i=1}^n\,\bigl( I(Y_i< U_i^{\rm T}\theta_m^*) - \tau_m \bigr)\,U_{imI},
\end{align}
and $U_{imI}$ is given in Section 3. Note that  $A_{nI}$ and $B_{nI}$ are the sub-vector(matrix) of $\mathbf A_n$ and $\mathbf B_n$ in $\eqref{0-anbn}$ corresponding to the index set $I$, respectively.

For $P_n(\gamma)$, define $p_n(\gamma) = P_n(\gamma)-P_n(\gamma^*)$.
Let $u_I$ be the subvector of $u$ corresponding to the subvector $\gamma_I$ in $\theta_I$. 
Then, we have
\begin{align*}
	p_n(\gamma_I,0) = p_n(\gamma_I^*+\nu_n u_I,0) 
	&= n\lambda_n\,\sum_{j\in J}\,\bigl(\sum_{m\in M_j}\,\omega_{mj}|\gamma_{mj}|\bigr)^{\frac{1}{2}}
		- n\lambda_n\,\sum_{j\in J}\,\bigl(\sum_{m\in M_j}\,\omega_{mj}|\gamma_{mj}^*|\bigr)^{\frac{1}{2}},
\end{align*}
which is $p_{1n}(\gamma)$ given in $\eqref{2-p1n}$. 
By $\eqref{2-tpn}$ in the proof of Theorem 1, we have 
$|p_{n}(\gamma_I,0)|= |p_n(\gamma_I^*+\nu_n u_I,0)|
	\leq n\lambda_nd_{nI}^\frac{1}{2}\nu_n^\frac{1}{2}s_n^\frac{3}{4}\|u\|_2^{\frac{1}{2}}$.
This, combined with $\eqref{4-qn}$, implies that
\begin{align}\label{4-ln-1}
	L_n(\theta_I^*+\nu_n u,0)-L_n(\theta_I^*,0)
	= \nu_n\mathbf{A}_{nI}^{\rm T}u+\frac{1}{2}\nu_n^2 u^{\rm T} \mathbf{B}_{nI} u + r_l(u),
\end{align}
where $\sup_{\|u\|_2\,\leq\,1}\,|r_l(u)| =O_p(n^{\frac{1}{2}}s_n^{\frac{3}{4}}{\nu_n}^{\frac{3}{2}})+O_p(ns_n^{\frac{1}{2}}\nu_n^3)
	+O(n\lambda_n\,d_{nI}^\frac{1}{2}\nu_n^\frac{1}{2}s_n^\frac{3}{4})$.

Define the ball $\Theta_{nI}=\{(\theta_I^{\rm T},0^{\rm T})^{\rm T}\,:\, \theta_I-\theta_I^*=\nu_n\,u\}$ with $\nu_n=C_0(n^{-1}p_n)^{\frac{1}{2}}$, where $u\in \mathbb{R}^{M+s_n}$ and $C_0$ is a positive constant.
Given that $\nu_n=C_0(n^{-1}p_n)^{\frac{1}{2}}$,  $\lambda_nd_{nI}^{\frac{1}{2}}=O(n^{-1}p_n^{\frac{1}{2}})$, and  $0<\alpha_0<\alpha_1<\frac{1}{6}$ we see that 
$n^{\frac{1}{2}}s_n^{\frac{3}{4}}{\nu_n}^{\frac{3}{2}} \leq Cn^{-\frac{1}{4}}s_n^{\frac{3}{4}}p_n^{\frac{3}{4}} = o_p(1)$,
$n\sqrt{s_n}\nu_n^3 \leq Cn^{-\frac{1}{2}}s_n^{\frac{1}{2}}p_n^{\frac{3}{2}} = o(n^{-\frac{1}{4}}s_n^{\frac{3}{4}}p_n^{\frac{3}{4}})= o_p(1)$
and $n\lambda_n\,d_{nI}^\frac{1}{2}\nu_n^\frac{1}{2}s_n^\frac{3}{4}\leq \lambda_n\,d_{nI}^\frac{1}{2}\,n^\frac{3}{4}s_n^\frac{3}{4}p_n^\frac{1}{4} = O(n^{-\frac{1}{4}}s_n^{\frac{3}{4}}p_n^{\frac{3}{4}})= o_p(1)$.
Hence, considering $\eqref{4-ln-1}$, for all $(\theta_I^{\rm T},0^{\rm T})^{\rm T}\in\Theta_{nI}$ we have
$L_n(\theta_I,0)-L_n(\theta_I^*,0)=\mathbf{A}_{nI}^{\rm T}(\theta_I-\theta_I^*)
	+\frac{1}{2}\,(\theta_I-\theta_I^*)^{\rm T}\,\mathbf{B}_{nI}\,(\theta_I-\theta_I^*) +  o_p(1)$,
which implies that for all $(\theta_I^{\rm T},0^{\rm T})^{\rm T}\in\Theta_{nI}$, 
\begin{eqnarray} \label{4-ln-eq}
    &&L_n(\theta_I,0)-L_n(\theta_I^*,0) + \frac{1}{2}\,\mathbf{A}_{nI}^{\rm T}\mathbf{B}_{nI}^{-1}\mathbf{A}_{nI}\nonumber\\ 
	&=& \frac{1}{2}\bigl(\mathbf{B}_{nI}^{-\frac{1}{2}}\,\mathbf{A}_{nI} + \mathbf{B}_{nI}^{\frac{1}{2}}(\theta_I-\theta_I^*)\bigr)^{\rm T}
	  \bigl(\mathbf{B}_{nI}^{-\frac{1}{2}}\,\mathbf{A}_{nI} + \mathbf{B}_{nI}^{\frac{1}{2}}(\theta_I-\theta_I^*)\bigr)
	  + o_p(1).  
\end{eqnarray}
By Theorem 1 and Theorem 2, we know that in probability a local minimizer $\hat{\theta}$ of  $L_n(\theta)-L_n(\theta^*)$ lies in the ball $\Theta_{nI}$,
which implies that $\hat{\theta}= (\hat{\theta}_I^{\rm T},0^{\rm T})^{\rm T} \in\Theta_{nI}$ in probability.
Hence, from $\eqref{4-ln-eq}$, we have
\begin{align}\label{4-11}
	\mathbf{B}_{nI}^{\frac{1}{2}}(\hat{\theta}_I-\theta_I^*) = -\mathbf{B}_{nI}^{-\frac{1}{2}}\,\mathbf{A}_{nI} + o_{p}(1)t,
\end{align}
where $t\in \mathbb{R}^{M+s_n}$ is an unit vector.
Since $t^{\rm T}\mathbf{B}_{nI}t\leq \lambda_{\max}(\mathbf{B}_{nI})\leq \lambda_{\max}(\mathbf{B}_n)\leq Cn$ by condition (L5), we have $\mathbf{B}_{nI}^{\frac{1}{2}}t=O_p(n^\frac{1}{2})t$.
Multiplying both sides of $\eqref{4-11}$ by a vector $b^{\rm T}\mathbf{B}_{nI}^{\frac{1}{2}}$, where $b \in \mathbb{R}^{M+s_n}$ is any unit vector, we obtain 
\begin{align}\label{4-13}
	b^{\rm T}\mathbf{B}_{nI}(\hat{\theta}_I-\theta_I^*) = -b^{\rm T}\mathbf{A}_{nI} + o_{p}(n^{1/2}).
\end{align}
Let $\xi_n=b^{\rm T}\mathbf{A}_{nI}$, and write $b=(b_1^{\rm T},\cdots,\,b_M^{\rm T})^{\rm T}$ where 
$b_m$ is the subvector of $b$ corresponding to the subvector $\theta_{mI}^*$ of $\theta_I^*$. 
By the definition of $\mathbf{A}_{nI}$ in $\eqref{4-an}$, we see that
\begin{align*}
	\xi_n &=\sum_{m=1}^M\,\sum_{i=1}^n\,\bigl(\psi_{mi}^*\,U_{imI}^{\rm T}b_m\bigr) = \sum_{i=1}^n\,\zeta_i,
\end{align*}
where $\zeta_i=\sum_{m=1}^M\,\psi_{mi}^*\,U_{imI}^{\rm T}b_m$ with $\psi_{mi}^* = I(Y_i <U_{imI}^{\rm T}\theta_{mI}^*)-\tau_m$. 
Clearly, $\{\zeta_i,\,i=1,\cdots, n\}$ is an independent sequence.
Next, we will verify that $\xi_n$ satisfies the Lindeberg's condition
\begin{align}\label{4-ldbg}
	\sigma_n^{-2}\,\sum_{i=1}^n\,E\bigl(\zeta_i^2I(|\zeta_i|\geq \sigma_n) \bigr) \to 0,
\end{align}
where $\sigma_n^2 = Var(\xi_n)$.

For $\zeta_i=\sum_{m=1}^M\,\psi_{mi}^*\,U_{imI}^{\rm T}b_m$, it is easy to see that $E(\zeta_i)=\sum_{m=1}^M\,E\bigl(E(\psi_{mi}^*|Z_i)\,U_{imI}^{\rm T}b_m\bigr)=0$ from the fact $E(\psi_{mi}^*|Z_i)=F(U_{imI}^{\rm T}\theta_{mI}^*|Z_i)-\tau_m=F(U_{i}^{\rm T}\theta_{m}^*|Z_i)-\tau_m=0$,
and that
\begin{align*}
	E(\zeta_i^2)
	&= E\bigl(\sum_{m=1}^M\,\psi_{mi}^*\,U_{imI}^{\rm T}b_m \bigr)^2
		= \sum_{m=1}^M\,\sum_{l=1}^M\,E\bigl(E(\psi_{mi}^*\psi_{li}^*|Z_i)\,b_m^{\rm T}(U_{imI}U_{ilI}^{\rm T})b_l\bigr)\notag\\
	&= \sum_{m=1}^M\,\sum_{l=1}^M\,E\bigl((\min(\tau_m,\tau_l) - \tau_m\tau_l)\,b_m^{\rm T}(U_{imI}U_{ilI}^{\rm T})b_l\bigr)
		= b^{\rm T}\mathbf{\Sigma}_nb,
\end{align*}
where $\mathbf{\Sigma}_n$ is given in $\eqref{4-sigma}$.
Hence, we obtain that $E(\xi_n)= \sum_{i=1}^n E(\zeta_i)=0$ and that, 
by independence $Var(\xi_n) = \sum_{i=1}^n E(\zeta_i^2)=nb^{\rm T}\mathbf{\Sigma}_nb$.
Under condition (L6), we obtain 
\begin{align}\label{4-sigma2}
	\sigma_n^2 = nb^{\rm T}\mathbf{\Sigma}_nb \geq n\lambda_{\min}(\mathbf{\Sigma}_n)b^{\rm T}b\geq Cn.
\end{align}
Note that $(U_{imI}^{\rm T}b_m)^2\leq \|U_{imI}\|_2^2\|b_m\|_2^2 \leq Cs_n\|b\|_2^2\leq Cs_n$ under condition (L2).
By Cauchy-Schwartz's inequality,  
$\zeta_i^2=(\sum_{m=1}^M\,\psi_{mi}^*\,U_{imI}^{\rm T}b_m)^2\leq M\sum_{m=1}^M\,(\psi_{mi}^*\,U_{imI}^{\rm T}b_m)^2\leq Cs_n\sum_{m=1}^M\,(\psi_{mi}^*)^2$,
which implies that $\zeta_i^2\leq Cs_n$ from the fact $|\psi_{mi}^*|\leq 1+\tau_m\leq 2$.
Hence, we have 
\begin{align*}
	\sum_{i=1}^n\,E\bigl(\zeta_i^2I(|\zeta_i|\geq \sigma_n) \bigr)
	&\leq Cs_n\sum_{i=1}^n\,E\bigl(I(\zeta_i^2\geq \sigma_n^2) \bigr)
		\leq Cs_n\sum_{i=1}^n\,P\bigl(Cs_n\sum_{m=1}^M\,(\psi_{mi}^*)^2 \geq \sigma_n^2\bigr)\notag\\
	&\leq Cs_n\sum_{i=1}^n\,E\bigl(Cs_n\sum_{m=1}^M\,(\psi_{mi}^*)^2\bigr)\,\sigma_n^{-2}
		\leq Cns_n^2\sigma_n^{-2}.
\end{align*}
Together with $\eqref{4-sigma2}$, this implies that, as $n\to \infty$,
\begin{align*}
	\sigma_n^{-2}\,\sum_{i=1}^n\,E\bigl(\zeta_i^2I(|\zeta_i|\geq \sigma_n) \bigr) \leq Cns_n^2\sigma_n^{-4} \leq C\frac{s_n^2}{n}\leq Cn^{-1+2\alpha_0} \to 0,
\end{align*}
which shows that Lindeberg's condition $\eqref{4-ldbg}$ holds. Hence,  
\begin{align*}
	(nb^{\rm T}\mathbf{\Sigma}_nb)^{-\frac{1}{2}}\,b^{\rm T}\mathbf{A}_{nI} 
	=\frac{b^{\rm T}\mathbf{A}_{nI}}{\sigma_n} = \frac{\xi_n-E(\xi_n)}{\sigma_n} \to N(0, 1).
\end{align*}
Together with $\eqref{4-13}$ and $\eqref{4-sigma2}$, this implies that
\begin{align*}
	(nb^{\rm T}\mathbf{\Sigma}_nb)^{-\frac{1}{2}}\,b^{\rm T}\mathbf{B}_{nI}\bigl(\hat{\theta}_I-\theta_I^* \bigr)
	=-(nb^{\rm T}\mathbf{\Sigma}_nb)^{-\frac{1}{2}}\,b^{\rm T}\mathbf{A}_{nI} +o_p(1) \to N(0, 1).
\end{align*}
This completes the proof of the theorem. \hfill$\Box$

%
%

%
%

%
%
%
%

\section*{References}

\bibliography{2015_10_06_mybib}

\begin{thebibliography}{23}
\expandafter\ifx\csname natexlab\endcsname\relax\def\natexlab#1{#1}\fi
\providecommand{\url}[1]{\texttt{#1}}
\providecommand{\href}[2]{#2}
\providecommand{\path}[1]{#1}
\providecommand{\DOIprefix}{doi:}
\providecommand{\ArXivprefix}{arXiv:}
\providecommand{\URLprefix}{URL: }
\providecommand{\Pubmedprefix}{pmid:}
\providecommand{\doi}[1]{\href{http://dx.doi.org/#1}{\path{#1}}}
\providecommand{\Pubmed}[1]{\href{pmid:#1}{\path{#1}}}
\providecommand{\bibinfo}[2]{#2}
\ifx\xfnm\relax \def\xfnm[#1]{\unskip,\space#1}\fi
\bibitem[{Avery et~al.(2011)Avery, He, North, Ambite, Boerwinkle, Fornage,
  Hindorff, Kooperberg, Meigs, Pankow et~al.}]{avery2011phenomics}
\bibinfo{author}{Avery, C.L.}, \bibinfo{author}{He, Q.},
  \bibinfo{author}{North, K.E.}, \bibinfo{author}{Ambite, J.L.},
  \bibinfo{author}{Boerwinkle, E.}, \bibinfo{author}{Fornage, M.},
  \bibinfo{author}{Hindorff, L.A.}, \bibinfo{author}{Kooperberg, C.},
  \bibinfo{author}{Meigs, J.B.}, \bibinfo{author}{Pankow, J.S.}, et~al.,
  \bibinfo{year}{2011}.
\newblock \bibinfo{title}{A phenomics-based strategy identifies loci on
  {APOC}1, {BRAP}, and {PLCG}1 associated with metabolic syndrome phenotype
  domains}.
\newblock \bibinfo{journal}{PLoS Genet} \bibinfo{volume}{7},
  \bibinfo{pages}{e1002322}.
\bibitem[{Fan and Li(2001)}]{fan2001variable}
\bibinfo{author}{Fan, J.}, \bibinfo{author}{Li, R.}, \bibinfo{year}{2001}.
\newblock \bibinfo{title}{Variable selection via nonconcave penalized
  likelihood and its oracle properties}.
\newblock \bibinfo{journal}{Journal of the American Statistical Association}
  \bibinfo{volume}{96}, \bibinfo{pages}{1348--1360}.
\bibitem[{Gasso et~al.(2009)Gasso, Rakotomamonjy and
  Canu}]{gasso2009recovering}
\bibinfo{author}{Gasso, G.}, \bibinfo{author}{Rakotomamonjy, A.},
  \bibinfo{author}{Canu, S.}, \bibinfo{year}{2009}.
\newblock \bibinfo{title}{Recovering sparse signals with a certain family of
  nonconvex penalties and {DC} programming}.
\newblock \bibinfo{journal}{IEEE Transactions on Signal Processing}
  \bibinfo{volume}{57}, \bibinfo{pages}{4686--4698}.
\bibitem[{Holland(2000)}]{holland2000glioblastoma}
\bibinfo{author}{Holland, E.C.}, \bibinfo{year}{2000}.
\newblock \bibinfo{title}{Glioblastoma multiforme: the terminator}.
\newblock \bibinfo{journal}{Proceedings of the National Academy of Sciences}
  \bibinfo{volume}{97}, \bibinfo{pages}{6242--6244}.
\bibitem[{Huang et~al.(2009)Huang, Ma, Xie and Zhang}]{huang2009group}
\bibinfo{author}{Huang, J.}, \bibinfo{author}{Ma, S.}, \bibinfo{author}{Xie,
  H.}, \bibinfo{author}{Zhang, C.H.}, \bibinfo{year}{2009}.
\newblock \bibinfo{title}{A group bridge approach for variable selection}.
\newblock \bibinfo{journal}{Biometrika} \bibinfo{volume}{96},
  \bibinfo{pages}{339--355}.
\bibitem[{Jiang et~al.(2014)Jiang, Bondell and Wang}]{jiang2014interquantile}
\bibinfo{author}{Jiang, L.}, \bibinfo{author}{Bondell, H.D.},
  \bibinfo{author}{Wang, H.J.}, \bibinfo{year}{2014}.
\newblock \bibinfo{title}{Interquantile shrinkage and variable selection in
  quantile regression}.
\newblock \bibinfo{journal}{Computational Statistics \& Data Analysis}
  \bibinfo{volume}{69}, \bibinfo{pages}{208--219}.
\bibitem[{Jiang et~al.(2013)Jiang, Wang and Bondell}]{jiang2013interquantile}
\bibinfo{author}{Jiang, L.}, \bibinfo{author}{Wang, H.J.},
  \bibinfo{author}{Bondell, H.D.}, \bibinfo{year}{2013}.
\newblock \bibinfo{title}{Interquantile shrinkage in regression models}.
\newblock \bibinfo{journal}{Journal of Computational and Graphical Statistics}
  \bibinfo{volume}{22}, \bibinfo{pages}{970--986}.
\bibitem[{Koenker(2004)}]{koenker2004quantile}
\bibinfo{author}{Koenker, R.}, \bibinfo{year}{2004}.
\newblock \bibinfo{title}{Quantile regression for longitudinal data}.
\newblock \bibinfo{journal}{Journal of Multivariate Analysis}
  \bibinfo{volume}{91}, \bibinfo{pages}{74--89}.
\bibitem[{Koenker(2015)}]{koenker2015quantile}
\bibinfo{author}{Koenker, R.}, \bibinfo{year}{2015}.
\newblock \bibinfo{title}{quantreg: Quantile Regression}.
\newblock \URLprefix \url{http://CRAN.R-project.org/package=quantreg}.
  \bibinfo{note}{{R} package version 5.11}.
\bibitem[{Koenker and Bassett(1978)}]{koenker1978regression}
\bibinfo{author}{Koenker, R.}, \bibinfo{author}{Bassett, J.G.},
  \bibinfo{year}{1978}.
\newblock \bibinfo{title}{Regression quantiles}.
\newblock \bibinfo{journal}{Econometrica} \bibinfo{volume}{46},
  \bibinfo{pages}{33--50}.
\bibitem[{Landmark-H{\o}yvik et~al.(2013)Landmark-H{\o}yvik, Dumeaux, Nebdal,
  Lund, Tost, Kamatani, Renault, B{\o}rresen-Dale, Kristensen and
  Edvardsen}]{landmark2013genome}
\bibinfo{author}{Landmark-H{\o}yvik, H.}, \bibinfo{author}{Dumeaux, V.},
  \bibinfo{author}{Nebdal, D.}, \bibinfo{author}{Lund, E.},
  \bibinfo{author}{Tost, J.}, \bibinfo{author}{Kamatani, Y.},
  \bibinfo{author}{Renault, V.}, \bibinfo{author}{B{\o}rresen-Dale, A.L.},
  \bibinfo{author}{Kristensen, V.}, \bibinfo{author}{Edvardsen, H.},
  \bibinfo{year}{2013}.
\newblock \bibinfo{title}{Genome-wide association study in breast cancer
  survivors reveals {SNP}s associated with gene expression of genes belonging
  to {MHC} class {I} and {II}}.
\newblock \bibinfo{journal}{Genomics} \bibinfo{volume}{102},
  \bibinfo{pages}{278--287}.
\bibitem[{Li and Zhu(2008)}]{li2008l1}
\bibinfo{author}{Li, Y.}, \bibinfo{author}{Zhu, J.}, \bibinfo{year}{2008}.
\newblock \bibinfo{title}{L$_1$-norm quantile regression}.
\newblock \bibinfo{journal}{Journal of Computational and Graphical Statistics}
  \bibinfo{volume}{17}, \bibinfo{pages}{163--185}.
\bibitem[{Mazumder et~al.(2011)Mazumder, Friedman and
  Hastie}]{mazumder2011sparsenet}
\bibinfo{author}{Mazumder, R.}, \bibinfo{author}{Friedman, J.H.},
  \bibinfo{author}{Hastie, T.}, \bibinfo{year}{2011}.
\newblock \bibinfo{title}{Sparsenet: Coordinate descent with nonconvex
  penalties}.
\newblock \bibinfo{journal}{Journal of the American Statistical Association}
  \bibinfo{volume}{106}, \bibinfo{pages}{1125--1138}.
\bibitem[{Peng et~al.(2014)Peng, Xu and Kutner}]{peng2014shrinkage}
\bibinfo{author}{Peng, L.}, \bibinfo{author}{Xu, J.}, \bibinfo{author}{Kutner,
  N.}, \bibinfo{year}{2014}.
\newblock \bibinfo{title}{Shrinkage estimation of varying covariate effects
  based on quantile regression}.
\newblock \bibinfo{journal}{Statistics and Computing} \bibinfo{volume}{24},
  \bibinfo{pages}{853--869}.
\bibitem[{Puputti et~al.(2006)Puputti, Tynninen, Sihto, Blom,
  M{\"a}enp{\"a}{\"a}, Isola, Paetau, Joensuu and
  Nupponen}]{puputti2006amplification}
\bibinfo{author}{Puputti, M.}, \bibinfo{author}{Tynninen, O.},
  \bibinfo{author}{Sihto, H.}, \bibinfo{author}{Blom, T.},
  \bibinfo{author}{M{\"a}enp{\"a}{\"a}, H.}, \bibinfo{author}{Isola, J.},
  \bibinfo{author}{Paetau, A.}, \bibinfo{author}{Joensuu, H.},
  \bibinfo{author}{Nupponen, N.N.}, \bibinfo{year}{2006}.
\newblock \bibinfo{title}{Amplification of {KIT}, {PDGFRA}, {VEGFR}2, and
  {EGFR} in gliomas}.
\newblock \bibinfo{journal}{Molecular Cancer Research} \bibinfo{volume}{4},
  \bibinfo{pages}{927--934}.
\bibitem[{Tibshirani(1996)}]{tibshirani1996regression}
\bibinfo{author}{Tibshirani, R.}, \bibinfo{year}{1996}.
\newblock \bibinfo{title}{Regression shrinkage and selection via the lasso}.
\newblock \bibinfo{journal}{Journal of the Royal Statistical Society: Series B}
  \bibinfo{volume}{58}, \bibinfo{pages}{267--288}.
\bibitem[{Tibshirani et~al.(2005)Tibshirani, Saunders, Rosset, Zhu and
  Knight}]{tibshirani2005sparsity}
\bibinfo{author}{Tibshirani, R.}, \bibinfo{author}{Saunders, M.},
  \bibinfo{author}{Rosset, S.}, \bibinfo{author}{Zhu, J.},
  \bibinfo{author}{Knight, K.}, \bibinfo{year}{2005}.
\newblock \bibinfo{title}{Sparsity and smoothness via the fused lasso}.
\newblock \bibinfo{journal}{Journal of the Royal Statistical Society: Series B}
  \bibinfo{volume}{67}, \bibinfo{pages}{91--108}.
\bibitem[{Wang et~al.(2012)Wang, Wu and Li}]{wang2012quantile}
\bibinfo{author}{Wang, L.}, \bibinfo{author}{Wu, Y.}, \bibinfo{author}{Li, R.},
  \bibinfo{year}{2012}.
\newblock \bibinfo{title}{Quantile regression for analyzing heterogeneity in
  ultra-high dimension}.
\newblock \bibinfo{journal}{Journal of the American Statistical Association}
  \bibinfo{volume}{107}, \bibinfo{pages}{214--222}.
\bibitem[{Wang et~al.(2009)Wang, Nan, Zhu and Zhu}]{wang2009hierarchically}
\bibinfo{author}{Wang, S.}, \bibinfo{author}{Nan, B.}, \bibinfo{author}{Zhu,
  N.}, \bibinfo{author}{Zhu, J.}, \bibinfo{year}{2009}.
\newblock \bibinfo{title}{Hierarchically penalized {C}ox regression with
  grouped variables}.
\newblock \bibinfo{journal}{Biometrika} \bibinfo{volume}{96},
  \bibinfo{pages}{307--322}.
\bibitem[{Wu and Liu(2009)}]{wu2009variable}
\bibinfo{author}{Wu, Y.}, \bibinfo{author}{Liu, Y.}, \bibinfo{year}{2009}.
\newblock \bibinfo{title}{Variable selection in quantile regression}.
\newblock \bibinfo{journal}{Statistica Sinica} \bibinfo{volume}{19},
  \bibinfo{pages}{801--817}.
\bibitem[{Zhao and Xiao(2014)}]{zhao2014efficient}
\bibinfo{author}{Zhao, Z.}, \bibinfo{author}{Xiao, Z.}, \bibinfo{year}{2014}.
\newblock \bibinfo{title}{Efficient regressions via optimally combining
  quantile information}.
\newblock \bibinfo{journal}{Econometric Theory} \bibinfo{volume}{30},
  \bibinfo{pages}{1272--1314}.
\bibitem[{Zou(2006)}]{zou2006adaptive}
\bibinfo{author}{Zou, H.}, \bibinfo{year}{2006}.
\newblock \bibinfo{title}{The adaptive lasso and its oracle properties}.
\newblock \bibinfo{journal}{Journal of the American Statistical Association}
  \bibinfo{volume}{101}, \bibinfo{pages}{1418--1429}.
\bibitem[{Zou and Yuan(2008)}]{zou2008regularized}
\bibinfo{author}{Zou, H.}, \bibinfo{author}{Yuan, M.}, \bibinfo{year}{2008}.
\newblock \bibinfo{title}{Regularized simultaneous model selection in multiple
  quantiles regression}.
\newblock \bibinfo{journal}{Computational Statistics \& Data Analysis}
  \bibinfo{volume}{52}, \bibinfo{pages}{5296--5304}.

\end{thebibliography}

\end{document}